
\font\tentitre=cmbx10 scaled \magstep2 
\font\tenmathtitre=cmbxti10 scaled\magstep2
\font\sevenmathtitre=cmbx7 scaled\magstep2
\newfam\mathtitrefam
\textfont\mathtitrefam=\tenmathtitre
\scriptfont\mathtitrefam=\sevenmathtitre
\def\bftitre{\fam\mathtitrefam\tentitre}

\font\tenbfmath=cmbxti10
\font\sevenbfmath=cmbx7
\newfam\bfmathfam
\textfont\bfmathfam=\tenbfmath\scriptfont\bfmathfam=\sevenbfmath
\def\bfmath{\fam\bfmathfam\tenbfmath}

\font\tenba=msbm10 scaled 900
\font\sevenba=msbm7 scaled 900
\font\fiveba=msbm5 scaled 900
\newfam\bafam
\textfont\bafam=\tenba\scriptfont\bafam=\sevenba
\scriptscriptfont\bafam=\fiveba
\def\ba{\fam\bafam\tenba}

\def\qed{\vbox{\hrule\hbox to 5pt{\vrule height 5pt\hfil\vrule}\hrule}}

\magnification=1200

{\nopagenumbers
\vbox{}
\vfill\vfill
\centerline{\bftitre
OPERATORS BETWEEN SUBSPACES} 
\vskip 1cm
\centerline{\bftitre AND QUOTIENTS OF ${\bftitre L^1}$}
\vfill
\centerline{\it G. GODEFROY - N.J. KALTON - D. LI}
\vfill\vfill\vfill\vfill\vfill
\eject
\vbox{}
\vfill
\centerline{\bftitre
OPERATORS BETWEEN SUBSPACES} 
\vskip 1cm
\centerline{\bftitre AND QUOTIENTS OF ${\bftitre L^1}$}
\vfill\vfill
\centerline{\it G. GODEFROY - N.J. KALTON\footnote{$^\dagger$} 
{\fiverm The second author was supported in part by NSF grant DMS-9870027}
- D. LI}
\vfill\vfill\vfill\vfill\vfill
\noindent{\it Abstract.} We provide an unified approach of results of L. Dor 
on the complementation of the range, and of D. Alspach on the nearness from 
isometries, of small into isomorphisms of $L^1$. We introduce the notion of 
small subspace of $L^1$ and show lifting theorems for operators between 
quotients of $L^1$ by small subspaces. We construct a subspace of $L^1$ 
which shows that extension of isometries from subspaces of $L^1$ to the 
whole space are no longer true for isomorphisms, and that nearly isometric 
isomorphisms from subspaces of $L^1$ into $L^1$ need not be near from any 
isometry.\par
\vskip 3mm
\noindent{\it R\'esum\'e.} On donne une m\'ethode qui permet d'obtenir en 
m\^eme temps le r\'esultat de L. Dor sur la compl\'ementation des 
sous-espaces de $L^1$ qui sont isomorphes \`a $L^1$ avec une petite 
constante d'isomorphisme, et le r\'esultat de D. Alspach disant qu'un 
isomorphisme de $L^1$ avec une petite constante d'isomorphisme est proche 
d'une isom\'etrie. Apr\`es avoir introduit la notion de ``petit'' 
sous-espace de $L^1$, on montre des th\'eor\`emes de rel\`evements pour les 
isomorphismes entre quotients de $L^1$ par des petits sous-espaces. On 
construit ensuite un sous-espace de $L^1$ qui montre que les th\'eor\`emes 
d'extension des isom\'etries entre sous-espaces de $L^1$ ne s'\'etendent pas 
au cas des isomorphismes, et que l'on peut avoir des isomorphismes de 
sous-espaces de $L^1$ \`a valeurs dans $L^1$ de constante arbitrairement 
proche de $1$, qui restent loin de toute isom\'etrie.\par
\vfill\vfill\vfill\vfill\vfill
\noindent{\bf Key-words.} {\it Integral representation of operators on $L^1$ 
-- Small into-isomorphisms of $L^1$ -- Complemented subspaces of $L^1$ -- 
Isometries of $L^1$ -- Small subspaces of $L^1$ -- Nicely placed subspaces 
of $L^1$ -- Semi-Riesz sets --  Strong Enflo operators -- Daugavet 
operators -- Lifting of operators -- $p$-stable random variables.}\par
\medskip
\noindent{\bf Mots-cl\'e.} {\it Repr\'esentation int\'egale des op\'erateurs 
de $L^1$ -- Isomorphismes de $L^1$ de petites constantes -- Sous-espaces 
compl\'ement\'es de $L^1$ -- Isom\'etries de $L^1$ -- Petits sous-espaces 
de $L^1$ -- Sous-espaces bien dispos\'es de $L^1$ -- Ensembles semi-Riesz 
-- Op\'erateurs d'Enflo forts -- Op\'erateurs de Daugavet -- 
Rel\`evements d'op\'erateurs -- Variables al\'eatoires $p$-stables.}\par
\medskip
\noindent {\bf AMS classification numbers}$\,:$ {\it Primary} 46A22 -- 
46B20 -- 46B25 -- {\it Secondary} 43A15 -- 43A46 -- 60E07.

\eject}

\pageno=1
\centerline{\bf OPERATORS BETWEEN SUBSPACES AND QUOTIENTS OF 
${\bfmath L^1}$}
\vskip 1cm
\centerline{\it G. GODEFROY - N.J. KALTON - D. LI}
\vskip 2cm
\noindent{\bf I. INTRODUCTION.}\par
\vskip 3mm
The Banach space $L^1$ is of fundamental importance for Fourier analysis, 
probability theory, and many other fields of pure and applied analysis. 
However, its Banach space structure is not yet fully understood. In 
particular, the study of operators from $L^1$ to $L^1$ lead to simply 
stated but apparently very hard problems$\,:$ for instance, it is not known 
whether a complemented subspace of $L^1$ is necessarily isomorphic to $L^1$ 
or $\ell^1$ (see [35]). A major tool for handling operators on $L^1(\Omega)$ 
is their representation through a ``matrix'' whose ``rows'' are indexed by 
the points of $\Omega$ and consist of measures on $\Omega$ ([18]; see also 
[9]). This representation will play a major role in the present work, 
together with duality arguments which will allow us to use methods and 
results from the isometric theory of Banach spaces. We will focus on 
several natural aspects of operator theory on $L^1$: perturbation of 
operators which are close to isometries, lifting of operators between 
quotient spaces, extensions of operators defined on linear subspaces. Such 
questions have of course been considered before and important theorems have 
been obtained (see e.g. [1], [7], [8], [15], [19], [21], [32], and 
references therein). Our work provides on one hand alternative simpler 
proofs to some known theorems, which are subsequently improved, and on the 
other hand new statements. We now turn to a detailed description of our 
results.\par
\vskip 3mm
Section II begins with a simple inequality, whose proof relies on Gaussian 
randomization, which provides a control of the atomic part of 
norm-increasing operators (Lemma II.1). This inequality leads to a unified 
approach to Dor's complementation theorem (Theorem II.3) and to Alspach's 
perturbation theorem (Theorem II.7). These two results follow from a 
Hahn-Banach argument applied to the atomic part of the relevant operator 
considered as taking values in a vector-valued $L^1$. Let us mention that 
if the limiting constant of $2^{1/2}$ obtained by L. Dor is not increased, 
our approach provides a quantitative improvement of Alspach's theorem. The 
new notion of ``small subspace'' is introduced in Section III (Definition 
III.1). Roughly speaking, a ``small subspace'' is a subspace of 
$L^1$ which is nowhere ``locally equal'' to $L^1$. Operators which have a 
non trivial ``diagonal'' in their ``matrix'' representation are denoted 
strong Enflo operators (Definition III.3). The link between these two 
notions is that operators which take their values in a small subspace are 
not strong Enflo$\,;$ in fact more is true (Proposition III.6). Small 
subspaces are studied within the class of translation invariant 
subspaces$\,;$ it is also shown that a closed direct sum of two small 
subspaces is small (Proposition III.9). Smallness of spaces is a crucial 
concept in section IV, which is devoted to the lifting of operators between 
quotient spaces$\,:$ the main result of this section is that if $X$ and $Y$ 
are small subspaces and the quotient spaces $L^1/X$ and $L^1/Y$ are close 
to each other in Banach-Mazur distance, then there is an isomorphism of 
$L^1$ which maps $X$ close to $Y$ with respect to the Hausdorff distance 
(Theorem IV.3). When one of the spaces $X$ or $Y$ has the lifting property 
(e.g. if one of  these spaces is complemented in its bidual) it follows 
that there is an invertible operator of $L^1$ which maps $X$ onto $Y$ 
(Corollary IV.8). This result can be seen as a substitute in $L^1$ of a 
theorem of Lindenstrauss and Rosenthal on $\ell^1$ ([24]$\,;$ see [25], 
Th.2.f.8). Note however that the quantitative behaviour of the liftings 
is actually better in the $L^1$ case (see Remark IV.9). Finally, section V 
provides a counterexample showing the ``sensibility over $\varepsilon$'' of 
the positive results obtained so far about isometries$\,:$ there is a 
subspace $X$ of $L^1$ such that, although isometries from this space into 
$L^1$ extend to isometries of the whole space $L^1$, there are linear maps 
$T$ from $X$ into $L^1$ which are arbitrarily close to isometries and such 
that no small perturbation of $T$ is actually an isometry from $X$ into 
$L^1\,;$ in particular, these ``near isometries'' do not extend to near 
isometries defined on $L^1$ (Theorem V.1). This space $X$ is arbitrarily 
close in the Banach-Mazur sense to weak-star closed subspaces of $\ell^1$, 
though it is far in the Hausdorff distance of any subspace of $L^1$ spanned 
by disjoint functions$\,;$ hence, although it has a strong ``$\,\ell^1\,$'' 
behaviour (for instance, its unit ball is compact locally convex in measure, 
see [12]), it cannot be decomposed into essentially disjoint finite 
dimensional parts (compare with [7], [33]). The construction of $X$ makes a 
crucial use of $p$-stable random variables with $p$ close to 1.\par
\vskip 5mm
Our notation is classical. All the Banach spaces we consider are 
indifferently real or complex, except in Section V, where we work only in 
real case. All the measure spaces are separable and purely nonatomic, and 
the measures are (unless other stated) probability measures$\,;$ hence the 
corresponding $L^1$-spaces are isometric to $L^1(0,1)$, but for convenience 
we use sometimes a more general probability space. When writing 
$L^1([0,1],m)$, we understand that $m$ is Lebesgue measure.\par\vskip 3mm 
The following representation of operators on $L^1$ ([18], Theorem 3.1) will 
be very useful$\,:$ for every operator $T\colon L^1\to L^1$, there are 
measures $\mu_x$, with $x\mapsto \mu_x$ weak$^\ast$-measurable, such 
that, for almost every $x$, 
\vskip -3mm 
$$Tf(x)=\int_0^1f(s)\,d\mu_x(s)=\sum_{n=1}^{+\infty} a_n(x)f(\sigma_n(x))+
\int_0^1f(s)\,d\nu_x(s)\leqno(R)$$\vskip -2mm\noindent
where $\displaystyle\smash{\sum_{n=1}^{+\infty}}a_n(x)\delta_{\sigma_n(x)}$ 
and $\nu_x$ are the atomic and continuous part of $\mu_x$.\par
\vskip 5mm
Section III and IV are essentially independent of Section II.
\par\vskip 5mm
\noindent{\bf Acknowledgement.} This work was initiated while the first 
author was Visiting Professor at the University of Missouri-Columbia during 
the academic year 1996--1997, and the third author was visiting this 
University in the fall of 1996. It is their pleasure to thank the 
Department of mathematics of U.M.C. for its hospitality and support, and 
all persons and institutions which made these stays possible.\par
\vfill\eject




\def\N{{\rm I \kern-.20em N}}
\def\ind{{\rm 1\kern-.30em I}}
\def\dis{\displaystyle}
\def\prob{\mathop{{\ba P}\kern 0pt}}

\magnification=1200

\noindent{\bf II. PROJECTIONS AND SMALL INTO-ISOMORPHISMS OF 
${\bfmath L^1}$.}
\vskip 3mm
In this section, we provide a unified approach to results of D. Alspach and 
L. Dor on isomorphisms from $L^1$ into $L^1$ with small bounds. A short 
proof of the special case of convolutors, which points towards the main 
lemma, will be given at the end of the section (see Remark II.8). The key 
result of this approach is the following lemma$\,:$\par
\vskip 5mm
\noindent{\bf Lemma II.1.} {\it Let $T\colon L^1([0,1],m)\to 
L^1([0,1],m)$ and write
$$Tf(s)=\sum_{n=1}^\infty a_n(s)f(\sigma_n(s))+
\int_0^1f(t)\,d\nu_s(t)\ .$$\par
Suppose that $\Vert Tf\Vert_1\geq\Vert f\Vert_1$ for all $f\in L^1.$ 
Then  for all $f\in L^1$   
$$\Vert f\Vert_1\leq\big\Vert (\sum_{k=1}^\infty \vert a_k(s)\vert^2
\vert f(\sigma_k(s))\vert^2)^{1/2}\big\Vert_1\ .$$}\par
\vskip 3mm
\noindent{\bf Proof.} Let $(\gamma_k)$ be a sequence of i.i.d. normal 
Gaussian (real or complex) random variables on a probability space 
$(\Omega,\prob)$.\par 
For any $l>0,$ we let 
$$D_k^l=\left[{k-1\over 2^l},{k\over 2^l}\right[\ ,\hskip 5mm 
1\leq k\leq 2^l,$$
be the dyadic intervals of level $l$.
\vskip 3mm
\noindent{\bf Fact.} {\it Let $\nu$ be a diffuse measure on $[0,1]$. 
Defining \vskip -2mm
$$\varphi_l=\sum_{k=1}^{2^l}\left\vert\nu(D_k^l)\right\vert^2\ ,$$
one has $\lim_l \varphi_l=0\,.$}\par
\vskip 3mm
Indeed, 
$$\varphi_l\leq\Vert \nu\Vert_1.\max_k
\big(\vert\nu(D_k^l)\vert\big)$$
since 
$(\sum\vert\lambda_k\vert^2)\leq(\sum\vert\lambda_k\vert).
(\max\vert\lambda_k\vert)$ for every sequence $(\lambda_k)$ of scalars. 
Since $\nu$ is diffuse, this $\dis{\max_k}$ tends to 0 when $l$ goes to 
infinity and the conclusion follows.\hfill$\qed$
\vskip 3mm
Let now $\rho$ be an arbitrary measure on $[0,1]$. Let us write 
$$\rho=\sum_{n=1}^\infty a_n\delta_{x_n}+\nu=\alpha+\nu$$
where $\alpha$ is atomic and $\nu$ is diffuse. For any $l$ and $N=2^l$, 
one has 
$$\sum_{k=1}^N \vert \rho(D_k^l)\vert^2=\sum_{k=1}^N\vert\alpha(D_k^l)+
\nu(D_k^l)\vert^2$$
and from the Fact, it follows easily that 
$$\lim_N\sum_{k=1}^N\vert\rho(D_k^l)\vert^2=
\sum_{n=1}^\infty\vert a_n\vert^2\ .\leqno(1)$$\par
We now denote $\dis{D_k^l=D_k},$ and for a given $N=2^l$, any 
$\omega\in\Omega$ and any $f\in L^1$, we define
$$f_\omega=\sum_{k=1}^N\gamma_k(\omega)\ind_{D_k}f\ .$$\par
With this notation, one has
$$\eqalign{\int_\Omega\Vert f_\omega\Vert_1\,d\prob(\omega)&=
\int\!\!\!\int_{[0,1]\times\Omega}\vert f_\omega(t)\vert\,dt\,d\prob(\omega)
=\int\!\!\!\int_{[0,1]\times\Omega}\vert\sum_{k=1}^N \gamma_k(\omega)
\ind_{D_k}(t)f(t)\vert\,dt\,d\prob(\omega)\cr
&=\int_\Omega\Big(\sum_{k=1}^N\vert \gamma_k(\omega)\vert
\int_{D_k}\vert f(t)\vert\,dt\Big)\,d\prob(\omega)=\Vert f\Vert_1\,,}$$
and it follows from the assumption on $T$ that for any $f$, one has
$$\Vert f\Vert_1\leq\int_\Omega\Vert T(f_\omega)\Vert_1\,d\prob(\omega)\,.$$
\par
We now use the representation $(R)$ of $T$
$$Tf(s)=\int_0^1f(t)\,d\mu_s(t)=\sum_{n=1}^\infty a_n(s)f(\sigma_n(s))+
\int_0^1f(t)\,d\nu_s(t)\ ,\leqno(R)$$
and for any given $f\in L^1$ we define for all $s$ a measure $\rho_s$ by 
$d\rho_s=f.d\mu_s\,.$ With this notation, one has 
$$T(f_\omega)(s)=\sum_{k=1}^N\gamma_k(\omega)\rho_s(D_k)$$
and thus by Fubini theorem
$$\int_\Omega\Vert T(f_\omega)\Vert_1\,d\prob(\omega)=
\int_0^1\!\!\int_\Omega\big\vert\sum_{k=1}^N
\gamma_k(\omega)\rho_s(D_k)\big\vert\,d\prob(\omega)\,ds=
\int_0^1\big(\sum_{k=1}^N\vert\rho_s(D_k)\vert^2\big)^{1/2}\,ds\,.$$\par
It follows that 
$$\Vert f\Vert_1\leq\int_0^1
\big(\sum_{k=1}^N\vert\rho_s(D_k)\vert^2\big)^{1/2}\,ds\,.$$\par
On the other hand, by applying $(1)$ to $\rho=\rho_s$, we have for almost 
every $s\,:$
$$\lim_N \big(\sum_{k=1}^N\vert\rho(D_k)\vert^2\big)^{1/2}
=\big(\sum_{n=1}^\infty\vert a_n(s)f(\sigma_n(s))\vert^2\big)^{1/2}\,.$$\par
Since $(R)$ implies that $\Vert\rho_s\Vert\leq 
\vert T\vert(\vert f\vert)(s)$, 
we may apply the dominated convergence theorem, and it follows then that
$$\Vert f\Vert_1\leq\Big\Vert\big(\sum_{k=1}^\infty\vert a_k(s)\vert^2
\vert f(\sigma_k(s))\vert^2\big)^{1/2}\Big\Vert_1\,,$$
which is Lemma 1.\hfill$\qed$\par
\vskip 5mm
\noindent{\bf Lemma II.2.} {\it Let $T\colon L^1([0,1],m)\to 
L^1([0,1],m)$ and $\alpha\geq1$ be such that 
$\alpha\Vert f\Vert_1\geq\Vert Tf\Vert_1\geq\Vert f\Vert_1$ for 
all $f\in L^1.$ Write 
$$Tf(s)=\sum_{n=1}^\infty a_n(s)f(\sigma_n(s))+
\int_0^1f(t)\,d\nu_s(t)\ .$$\par
Then for all $f\in L^1$
$$\big\Vert\sum_{n=1}^\infty\vert a_n(s)\vert\,\vert f(\sigma_n(s))\vert\,
\big\Vert_1\leq\alpha\Vert f\Vert_1\leqno(\ast)$$
$$\alpha^{-1}\Vert f\Vert_1\leq\Vert\max_n \vert a_n(s)\vert\,
\vert f(\sigma_n(s))\vert\,\Vert_1 \ . 
\leqno(\ast\ast)$$}\vskip 3mm
\noindent{\bf Proof.} It follows e.g. from $(R)$ that
$\Vert\,\vert T\vert\,\Vert=\Vert T\Vert,$ hence $\Vert\,\vert T\vert (f)
\Vert_1\leq \alpha\Vert f\Vert_1,$ which shows $(\ast).$\par
For $(\ast\ast)$, we observe that for any $s\in [0,1]$
$$\big(\sum_{k=1}^\infty \vert a_k(s)\vert^2\vert f(\sigma_k(s))\vert
^2\big)^{1/2}\leq\big(\sum_{k=1}^\infty\vert a_k(s)\vert\,\vert 
f(\sigma_k(s)\vert\big)^{1/2}\big(\max_k \vert a_k(s)\vert\,
\vert f(\sigma_k(s))\vert\big)^{1/2}\ ,$$
hence by Cauchy-Schwarz inequality
$$\eqalign{\Big\Vert \big(\sum_{k=1}^\infty \vert a_k(s)\vert^2\vert 
f(\sigma_k(s))\vert^2\big)&^{1/2}\Big\Vert_1\cr
\leq \Big\Vert
&\big(\sum_{k=1}^\infty\vert a_k(s)\vert\,\vert 
f(\sigma_k(s)\vert\big)\Big\Vert_1^{1/2}.\,\Vert\max_k \vert a_k(s)\vert\,
\vert f(\sigma_k(s))\vert\,\Vert_1^{1/2}.}$$\par
It follows now from $(\ast)$ and Lemma 1 that
$$\Vert f\Vert_1\leq\sqrt\alpha \Vert f\Vert_1^{1/2}
\Vert\max_k\vert a_k(s)\vert\,\vert f(\sigma_k(s)\vert\,\Vert_1^{1/2}$$
which shows $(\ast\ast)\ .$\hfill$\qed$\par
\vskip 5mm
Lemma II.2 allows us to show the following theorem due to L. Dor ([7], 
Corol. 3.3) through a Hahn-Banach argument, distinct however from the 
original argument ([7]). Note that we do not have the better constant 1.6 
obtained by L. Dor in the real case.\par
\vskip 5mm
\noindent{\bf Theorem II.3.} {\it Let $T\colon L^1\to L^1$ be an 
isomorphism into $L^1$ such that $\Vert T\Vert.\Vert T^{-1}\Vert<\sqrt2,$ 
with $T^{-1}\colon T(L^1)\to L^1.$ Then $T(L^1)$ is a complemented 
subspace of $L^1.$}\par
\vskip 3mm
\noindent{\bf Proof.} We may and do assume that $T$ satisfies the 
assumptions of Lemma II.2 with $\alpha<2$.\par
Define $\Psi\colon L^1\to L^1(c_0)$ by
$$\Psi(f)=\left(\vert a_n(s)\vert\,f(\sigma_n(s))\right)_{n\geq1}\ .$$
Note that, by $(\ast),$ $\Psi$ actually takes its values in $L^1(\ell_1).$ 
By $(\ast\ast),$ we have
$$\Vert \Psi(f)\Vert_{L^1(c_0)}\geq \alpha^{-1}\Vert f\Vert_1
\hskip 2cm \forall f\in L^1\ .$$
Hence, there exists $F\in L^1(c_0)^\ast=L^\infty(\ell_1)$  with 
$F\geq 0,$ $\Vert F\Vert\leq\alpha,$ such that
$$\int_0^1 f(s)\,ds=F[\Psi(f)]\hskip 2cm \forall f\in L^1\,.$$
That means there exists a sequence $(b_n)_{n\geq1}$ of positive functions 
such that
$$\Vert \sum_{n=1}^\infty b_n(s)\Vert_\infty\leq\alpha\ ,$$
and
$$\int_0^1f(s)\,ds=\int_0^1\sum_{n=1}^\infty b_n(s)\,\vert a_n(s)\vert
\,f(\sigma_n(s))\,ds\hskip 1cm\forall f\in L^1\ .$$
Choose disjoint measurable subsets $E_n,$ $n\geq1,$ of $[0,1]\times[0,1]$ 
such that
$$m(\{x\,;\,(s,x)\in E_n\})={b_n(s)\over \alpha}\hskip 1cm \forall 
s\in[0,1],\hskip 3mm \forall n\geq1\,.$$
We define $U\colon L^1([0,1])\to L^1([0,1]^2)$ by
$$Uf(s,x)=\sum_{n=1}^\infty a_n(s)\,f(\sigma_n(s))\,\ind_{E_n}(s,x)\ .$$
For any $f\in L^1,$ one has
$$\Vert Uf\Vert_1=\int_0^1\sum_{n=1}^\infty{b_n(s)\over\alpha}\vert 
a_n(s)\vert\,\vert f(\sigma_n(s))\vert\,ds={1\over\alpha}\Vert f\Vert
_1\ .$$
Hence $U(L^1)$ is isometric to $L^1$ and it follows (see [22], Chap. 6, 
\S 17, Theorem 3) that $U(L^1)$ is $1$-complemented in $L^1([0,1]^2).$\par
\vskip 3mm
Define $\overline T\colon L^1([0,1])\to L^1([0,1]^2)$ by
$$\overline Tf(s,x)=Tf(s)=\sum_{n=1}^\infty a_n(s)\,f(\sigma_n(s))+\int_0^1 
f(u)\,d\nu_s(u)\ .$$\par
If we denote
$$\eqalign{Uf(s,x)&=\int_0^1f(u)\,d\mu_{(s,x)}^U(u)\cr
\overline Tf(s,x)&=\int_0^1f(u)\,d\mu_{(s,x)}^{\overline T}(u)\ ,}$$
then clearly
$$\mu_{(s,x)}^U=a_k(s)\delta_{\sigma_k(s)}$$
for the unique $k\geq1$ such that $(s,x)\in E_k,$ and
$$\mu_{(s,x)}^{\overline T}=\sum_{n=1}^\infty a_n(s)\delta_{\sigma_n(s)}
+\nu_s$$
for all $(s,x).$ Hence, we have
$$\mu_{(s,x)}^{\overline T-U}=\sum_{n\not=k} a_n(s)\delta_{\sigma_n(s)}
+\nu_s$$
and it easily follows that
$$\vert \overline T\vert-\vert U\vert=\vert\overline T-U\vert\geq0\ .$$\par
Note that for all $f\in L^1:$
$$\Vert\,\vert U\vert (f)\Vert_1=\Vert Uf\Vert_1={1\over\alpha}
\Vert f\Vert_1\,.$$\par
On the other hand,
$$\Vert\,\vert\overline T\vert (f)\Vert_1=\Vert\,\vert T\vert(f)\Vert_1
\leq\alpha\Vert f\Vert_1\,.$$\par
Hence, for any $f\geq0,$
$$\eqalign{\Vert\,\vert\overline T-U\vert(f)\Vert_1&=\Vert\,\vert
\overline T\vert(f)-\vert U\vert(f)\Vert_1=\Vert\,\vert\overline T\vert(f)
\Vert_1-\Vert\,\vert U\vert(f)\Vert_1\cr
&\leq(\alpha-\alpha^{-1})\Vert f\Vert_1=(\alpha^2-1)
\Vert Uf\Vert_1\,.}$$\par
Observe now that $U$ maps disjoint functions to disjoint functions, 
hence if $f=f^+-f^-$
$$\Vert Uf\Vert_1=\Vert Uf^+\Vert_1+\Vert Uf^-\Vert_1\,,$$
and thus for any $g\in L^1\,:$
$$\Vert (\overline T-U)g\Vert_1\leq(\alpha^2-1)\Vert Ug\Vert_1=\beta 
\Vert Ug\Vert_1\,.\leqno(2)$$\par
Since $\alpha<\sqrt2,$ we have $\beta=\alpha^2-1<1.$ Let 
$\Pi:L^1([0,1]^2)\to U(L^1)$ be a projection with $\Vert\Pi\Vert 
=1.$ It follows from $(2)$ that $\dis{\Pi_{\vert \overline T
(L^1)}}$ is an isomorphism from $\overline T(L^1)$ onto $U(L^1)\,.$
\par
Let $S=\Pi^{-1}\colon U(L^1)\to \overline T(L^1).$ It is easily seen that 
$S\Pi$ defines a projection from $L^1([0,1]^2)$ onto $\overline T(L^1).$ 
Restricting $S\Pi$ to the functions which are independent of $x$ gives a 
projection from $L^1$ to $T(L^1).$\hfill$\qed$\par
\vskip 5mm
Our goal is now to provide a similar approach to Alspach's result on near 
isometries ([1]). We first prove$\,:$\par
\vskip 5mm
\noindent{\bf Theorem II.4.} {\it Let $T\colon L^1\to L^1$ be such that
$$\alpha\Vert f\Vert_1\geq\Vert Tf\Vert_1\geq\Vert f\Vert_1\hskip 1cm
\forall f\in L^1\,.$$
Then for each $\varepsilon>0,$ there exists a 
$(\tau_m$-\/$\tau_m)$-continuous operator $S_\varepsilon\colon L^1\to L^1$ 
such that
$$\Vert T-S_\varepsilon\Vert\leq(\alpha-1)+\varepsilon\ .$$}\par
\vskip 3mm
\noindent{\bf Proof.} We consider again
$$\Psi(f)=\left(\vert a_n(s)\vert\,f(\sigma_n(s))\right)_{n\geq1}\,,$$
but this time as an operator from $L^1$ to $L^1(\ell_2).$ By Lemma 1, 
we have for all $f$
$$\Vert f\Vert_1\leq\Vert\Psi(f)\Vert_{L^1(\ell_2)}\ .$$
\par
Since $L^1(\ell_2)^\ast=L^\infty(\ell_2),$ there is $(b_n)_{n\geq1}$ with
$b_n\geq0$ such that
$$\left\{\eqalign{&\Vert(\sum b_n^2)^{1/2}\Vert_\infty\leq1\ ;\cr
&\int_0^1\sum b_n(s)\,\vert a_n(s)\vert\,f(\sigma_n(s))\,ds=\int_0^1 f\,ds}
\right.$$
for all $f.$\par
We may and do redefine the $a_n,b_n,\sigma_n$'s in such a way that 
$b_1\geq b_2\geq\cdots\,;$ hence $0\leq b_n(s)
\leq 1/\sqrt n$ for all $s$ and $n.$ Then
$$\int_0^1\vert\sum_{n=N+1}^\infty b_n(s)a_n(s)f(\sigma_n(s))\vert\,ds\leq
{1\over\sqrt N}\int_0^1\sum_{n=N+1}^\infty 
\vert a_n(s)\vert\,\vert f(\sigma_n(s))\vert\,ds\leq{1\over\sqrt N}
\Vert T\Vert.\Vert f\Vert_1$$
thus
$$\int_0^1\sum_{n=1}^N b_n(s)\vert a_n(s)\vert\,\vert f(\sigma_n(s))\vert\,
ds\geq\Vert f\Vert_1-{1\over\sqrt N}\Vert T\Vert.\Vert f\Vert_1\ .$$\par
We let
$$S_Nf(s)=\sum_{n=1}^N a_n(s)b_n(s)f(\sigma_n(s))\ .$$
For any $s$
$$\mu_s^{T-S_N}=\sum_{n=1}^N a_n(s)(1-b_n(s))\delta_{\sigma_n(s)}+\nu_s
+\sum_{n=N+1}^\infty a_n(s)\delta_{\sigma_n(s)}\ ;$$
therefore $\vert T-S_N\vert=\vert T\vert-\vert S_N\vert,$ and for any 
$f\geq0,$ one has
$$\Vert\,\vert T-S_N\vert(f)\Vert_1=
\Vert\,\vert T\vert(f)-\vert S_N\vert(f)\Vert_1
=\Vert\,\vert T\vert(f)\Vert_1-\Vert\,\vert S_N\vert(f)\Vert_1\ .$$
Since
$$\eqalign{&\Vert\,\vert T\vert(f)\Vert_1\leq\alpha\Vert f\Vert_1\cr
&\Vert\,\vert S_N\vert(f)\Vert_1\geq
\Vert f\Vert_1-{1\over\sqrt N}\Vert T\Vert.\Vert f\Vert_1}$$
one has for all $f\geq0$
$$\Vert\,\vert T-S_N\vert(f)\Vert_1\leq(\alpha-1)\Vert f\Vert_1+
{1\over\sqrt N}\Vert T\Vert.\Vert f\Vert_1$$
and thus
$$\Vert T-S_N\Vert=\Vert\,\vert T-S_N\vert\,\Vert\leq(\alpha-1)+
{\Vert T\Vert\over\sqrt N}\,.$$\par
The conclusion follows since $N$ is arbitrary and $S_N$ is 
$(\tau_m$-\/$\tau_m)$-continuous.\hfill$\qed$\par 
\vskip 5mm
\noindent{\bf Proposition II.5.} {\it Let $S\colon L^1\to L^1$ be an 
operator such that
$$Sf(s)=\sum_{k=1}^\infty a_k(s)f(\sigma_k(s))$$
and
$$\alpha\Vert f\Vert_1\geq\Vert Sf\Vert_1\geq\Vert f\Vert_1$$
for every $f\in L^1.$\par
Then there exists an operator $U\colon L^1\to L^1$ of the form
$$Uf(s)=c(s)f(\sigma(s))$$
such that 
$$\Vert S-U\Vert\leq 2(\alpha-{1\over\alpha})\ .$$}\par
\vskip 3mm
\noindent{\bf Proof.} We consider as before 
$\dis{\Psi(f)=\left(\vert a_n(s)\vert\,f(\sigma_n(s))\right)_{n\geq1}\,,}$ 
as an operator into $L^1(c_0).$ Since by $(\ast\ast)$
$$\alpha^{-1}\Vert \Psi(f)\Vert_{L^1(c_0)}\leq\Vert f\Vert_{L^1},$$ 
there exists $(b_k)_k$ with $b_k\geq0\,,$ $\sum_k b_k\leq1\,,$ and
$$\int_0^1\sum_{k=1}^\infty\vert a_k(s)\vert\, b_k(s) 
f(\sigma_k(s))\,ds={1\over\alpha}\int_0^1f(s)\,ds\ .$$
We may and do assume that $b_1\geq b_2\geq \cdots.$ Therefore, for all 
$k\geq2,$
$$1-b_k\geq 1-b_1\geq b_k\ .$$\par
Since $\alpha\Vert f\Vert_1\geq\Vert\,\vert S\vert(f)\Vert_1\,,$ one has,
for all $f\geq0,$
$$\eqalign{(\alpha-{1\over \alpha})\Vert f\Vert_1&\geq
\int_0^1\sum_{k=1}^\infty \vert a_k(s)\vert(1-b_k(s))f(\sigma_k(s))\,ds\cr
&\geq
\int_0^1\sum_{k=2}^\infty \vert a_k(s)\vert\,b_k(s)f(\sigma_k(s))\,ds\,.}$$
\par
Hence if we define
$$Vf(s)=\sum_{k=2}^\infty a_k(s)b_k(s)f(\sigma_k(s))$$
and
$$\eqalign{Wf(s)&=\sum_{k=1}^\infty a_k(s)b_k(s)f(\sigma_k(s))\cr
&=a_1(s)b_1(s)f(\sigma_1(s))+Vf(s)\ ,}$$
one has for all $f\geq0$
$$\Vert\,\vert S-W\vert(f)\Vert_1=\Vert\,\vert S\vert(f)-\vert W\vert(f)
\Vert_1\leq (\alpha-{1\over\alpha})\Vert f\Vert_1\ ,$$
and also
$$\Vert\,\vert V\vert(f)\Vert_1\leq (\alpha-{1\over\alpha})
\Vert f\Vert_1\ ;$$
thus 
$$\Vert S-W\Vert\leq (\alpha-{1\over\alpha})\hskip 5mm{\rm and}
\hskip 5mm\Vert V\Vert\leq (\alpha-{1\over\alpha}).$$\par
We let now$\,:$ $U=W-V.$\par
We have$\,:$ 
$$\Vert S-U\Vert=\Vert S-W+V\Vert\leq \Vert S-W\Vert+
\Vert V\Vert\leq 2(\alpha-{1\over\alpha})\ .\eqno\qed$$
\par\vskip 5mm\goodbreak
\noindent{\bf Proposition II.6.} {\it Let $U\colon L^1\to L^1$ be of the 
form
$$Uf(s)=c(s)f(\sigma(s))$$
and assume that for some $\alpha\geq1,$ we have for all $f\in L^1$
$$\alpha\Vert f\Vert_1\geq \Vert U(f)\Vert_1\geq \Vert f\Vert_1\ .$$
\par
Then, there exists an isometry $J$ from $L^1$ to $L^1$ such that$\,$  
$\Vert U-J\Vert\leq(\alpha-1)\ .$}\par
\vskip 3mm
\noindent{\bf Proof.} By the Radon-Nikodym theorem, there exists $w$ such 
that for all $f\in L^1$
$$\int_0^1\vert U\vert(f)(s)\,ds=\int_0^1w(s)f(s)\,ds\ .$$
\par
For all $f\geq0,$ one has
$$\alpha\int_0^1 f(s)\,ds\geq\int_0^1 w(s)f(s)\,ds\geq\int_0^1f(s)\,ds\ .$$
It follows that $1\leq w(s)\leq\alpha$ for all $s.$ We define
$$Jf(s)={c(s)\over w(\sigma(s))}f(\sigma(s))$$
and compute
$$\eqalign{\Vert Jf\Vert_1&=\int_0^1\vert Jf(s)\vert\,ds=
\int_0^1 {\vert c(s)\vert\over w(\sigma(s))}\vert f(\sigma(s))\vert\,ds\cr
&=\int_0^1\vert c(s)\vert\left({\vert f\vert\over w}\right)
(\sigma(s))\,ds\cr
&=\int_0^1 \vert U\vert\left({\vert f\vert\over w}\right)(s)\,ds
=\int_0^1 w(s)\left({\vert f\vert\over w}\right)(s)\,ds=\Vert f\Vert_1\ .}$$
\par
Hence $J$ is an isometry.\par
Moreover, for any $s,$ we have
$$\mu_s^{\vert U-J\vert}=\vert c(s)\vert\left(1-{1\over w(\sigma(s))}\right)
\delta_{\sigma(s)}\ ,$$
and since
$$0\leq 1-{1\over w(s)}\leq 1-{1\over \alpha}={\alpha-1\over\alpha}\,,$$
we have
$$\Vert U-J\Vert=\Vert\,\vert U-J\vert\,\Vert\leq{\alpha-1\over\alpha}
\Vert U\Vert \leq (\alpha-1)\ .\eqno\qed$$
\par\vskip 5mm\vfill
From Theorem II.4, Proposition II.5 and Proposition II.6 we deduce the 
following quantitative improvement of Alspach's theorem ([1]).\par
\vskip 5mm
\noindent{\bf Theorem II.7.} {\it There is a function $\varphi(\alpha)$ 
with $\dis{\lim_{\alpha\to1_+}\varphi(\alpha)=0}$ and such that if 
$T\colon L^1\to L^1$ satisties for some $\alpha\geq1$
$$\alpha\Vert f\Vert_1\geq \Vert Tf\Vert_1\geq \Vert f\Vert_1\hskip 1cm
\forall f\in L^1\,,$$
then there is an isometry $J\colon L^1\to L^1$ such that for all 
$f\in L^1$
$$\Vert T-J\Vert\leq \varphi(\alpha)\ .$$\par
Moreover, for $1\leq \alpha\leq 1.08$ we can choose 
$\varphi(\alpha)=13(\alpha-1).$}\par
\vskip 3mm
\noindent{\bf Proof.} The first part clearly follows from Theorem II.4, 
Proposition II.5 and Proposition II.6.\par {\openup 1mm
For the second part, apply Theorem II.4 to find an $S$ so that 
$\Vert T-S\Vert\leq {17\over 16}(\alpha-1)=\delta_1<1\,.$ Then 
$S'=\dis{1\over 1-\delta_1}S$ verify $\Vert f\Vert_1\leq\Vert S'f\Vert_1
\leq \alpha_1\Vert f\Vert_1,$ with 
$\alpha_1=\dis{\delta_1+\alpha\over 1-\delta_1}\leq1.28\,,$
and $\Vert T-S'\Vert\leq \dis{\delta_1\over 1-\delta_1}(1+\alpha)\leq 
2.5(\alpha-1).$ Now, Proposition II.5 gives an operator $U$ so that 
$\Vert S'-U\Vert\leq 2(\alpha_1-\alpha_1^{-1})=\delta_2<1\,;$ 
the operator 
$U'=\dis{1\over 1-\delta_2}U$ verify 
$\Vert f\Vert_1\leq\Vert U'f\Vert_1\leq\alpha_2\Vert f\Vert_1$, with 
$\alpha_2=\dis{\delta_2+\alpha_1\over 1-\delta_2}$ and we have 
$\Vert S'-U'\Vert \leq \dis{\delta_2\over 1-\delta_2}(1+\alpha_1) 
\leq 3.4(\alpha-1)\,.$\par\vskip 1mm
Finally, Proposition II.6 gives a $J$ so that 
$\Vert U'-J\Vert\leq \alpha_2-1\leq 6.8(\alpha-1).$\hfill
$\qed$}
\vskip 5mm
\noindent{\bf Remark II.8.} We give a proof of Theorem II.7 in the special 
case where $T=C_\mu$ is the convolution by a measure $\mu$ on 
$L^1({\ba T})$.\par
Suppose that $\Vert \mu\Vert=1$ and that $\Vert f\ast\mu\Vert_1\geq (1-
\varepsilon)\Vert f\Vert_1$. Write the atomic part of $\mu$ as 
$\dis{\sum_{k=1}^{+\infty}a_k\delta_{x_k}}$. We have 
$$\sum_{k=1}^{+\infty}\vert a_k\vert\leq \Vert\mu\Vert=1\,.\leqno(3)$$
On the other hand, Wiener's theorem ([20], p. 42$\,;$ [13], p. 415)$\,:$
$$\lim_{N\to +\infty}{1\over 2N+1}\sum_{n=-N}^N
\vert \widehat\mu(n)\vert^2=\sum_{k=1}^{+\infty}\vert a_k\vert^2$$ 
says that 
$$(1-\varepsilon)\leq\sum_{k=1}^{+\infty}\vert a_k\vert^2\leqno(4)$$
since 
$$\vert \widehat\mu(n)\vert=\Vert\mu\ast e_n\Vert_1\geq (1-\varepsilon)\,,$$ 
where $e_n(t)=e^{2\pi int}.$ It follows from $(3)$ and $(4)$ that there 
exists an index $K$ such that $\vert a_K\vert\geq (1-2\varepsilon)\,;$ then 
the distance between $C_\mu$ and the translation by $x_K$ of less than 
$4\varepsilon,$ and this conclude the proof.\hfill$\qed$\par
There is an interesting link between Wiener's theorem and the equation $(1)$ 
in the proof of Lemma II.1. Indeed, let $\mu$ be a measure on the Cantor 
group $G=\{-1,1\}^\N,$ and let $(W_k)$ and $(D_k)$ be the Walsh functions 
and the elementary open sets with index $k\in \{-1,1\}^{[\N]}.$ With this 
notation, one has by Parseval's formula\vskip -2mm
$$2^{-l}\sum_{k\in2^l}\big(\int W_k\,d\mu\big)^2=\sum_{k\in 2^l}\mu(D_k)^2$$
\vskip-1mm\noindent
for any $l\geq1,$ and thus equation $(1)$ appears, through the natural 
measure preserving isomorphism between $G$ and the unit interval, as the 
exact analogue for the Cantor group of Wiener's theorem. It would be nice 
to state a generalization of Wiener's theorem to arbitrary compact abelian 
groups.
\vskip 1cm





\def\N{{\rm I \kern-.25em N}}
\def\ind{{\rm 1\kern-.30em I}}
\def\dis{\displaystyle}
\def\dlowlim#1#2{\mathop{\rm #1}\limits_{#2}}
\def\ksubset{\mathop{\raise -1mm
\vbox{\hbox{$\subset$}\baselineskip=1mm\hbox{$\sim$}}}} 

\magnification=1200

\noindent{\bf III. SMALL SUBSPACES OF ${\bfmath L^1}$.}\par
\vskip 3mm
In this section, we introduce and study the notion of {\it small subspace} 
of $L^1$ which will be essential in the next section.\par\vskip 1mm
In the sequel the notation $\,:$ $A\ksubset D$ will be used to say that 
there is a real number $k>0$ such that $kA\subseteq D.$\par
If $X$ is a subspace of $L^1,$ we denote by $C_X=\overline{B_X}^{\tau_m}$ 
the closure of its unit ball for the topology $\tau_m$ of convergence in 
measure. Recall that Bukhvalov and Lozanovski ([4]) showed that 
$P\big(B_{X^{\perp\perp}}\big)=C_X,$ where $P$ is the natural 
projection from $L^{1\ast\ast}$ to $L^1.$\par
We also denote, for $A\subseteq \Omega$
$$L^1(A)=\{f\in L^1\,;\ \ind_{\Omega\setminus A}.f=0\hskip 3mm 
{\rm a.e.}\}\ .$$
\vskip 3mm 
\noindent{\bf Definition III.1.} {\it A subspace $X$ of $L^1(\Omega,m)$ is 
said to be {\sl small} if there is no $A\subseteq\Omega$ of positive 
measure such that $B_{L^1(A)}\ksubset \ind_A.B_X\ .$}\par\vskip 3mm
In other words, the projection $f\mapsto f.\ind_A$ never maps $X$ onto 
$L^1(A).$\par
It is easy to see (by lifting finite trees for instance) that reflexive 
subspaces of $L^1$ are small. However smallness is not an intrinsic notion, 
but is related to the position of $X$ in $L^1$. Indeed, let 
$Q\colon\ell_1\to L^1$ be a quotient map and $J\colon\ell_1\to L^1$ be an 
isometric embedding. Then, the range of 
$(J,Q)\colon\ell_1\to L^1\times L^1$ is isomorphic to $\ell_1$ but is not 
small. On the other hand, in $L^1([0,1]\times[0,1])$ the subspace of the 
functions which do not depend of the second coordinate is small, although 
it is isometric to $L^1.$\par  
We are now going to give some examples and properties of small subspaces.
\par 
It will be useful to know that small subspaces actually satisfy a formally 
stronger property.\par
\vskip 3mm
\noindent{\bf Proposition III.2.} {\it If $X$ is a small subspace of $L^1,$ 
there is no $E\subseteq\Omega$ of positive measure such that 
$B_{L^1(E)}\ksubset \ind_E.C_X\ .$}\par
\vskip 3mm
\noindent{\bf Proof.} Assume that $B_{L^1(E)}\subseteq \ind_E.(kC_X)\,,$ 
and let $\{x_n\}$ be dense in $B_{L^1(E)}.$ There are $f_n\in kB_X,$ and 
$E_n\subseteq E$ with $m(E_n)\leq 4^{-n}m(E)$ such that 
$\vert x_n-f_n\vert\leq 1/n$ on $(E\setminus E_n).$ If now 
$A=E\setminus\big(\cup_{n\geq1}E_n\big),$ we have 
$B_{L^1(A)}\subseteq \ind_A.(kB_X)\,,$ and so $X$ is not small.\hfill
$\qed$\par
\vskip 5mm
The following notion is closely related to smallness. It should be however 
noted that, except in the proof of Proposition III.9, only the property 
stated in Proposition III.4 will actually be used.\par
\vskip 5mm
\noindent{\bf Definition III.3.} {\it Let $T$ be an operator on $L^1,$ 
written as
$$Tf(x)=\int_0^1 f(s)\,d\mu_x(s)\ .$$
$T$ will be called a {\sl strong Enflo operator} if there exists a set 
of positive measure on which $\mu_s(\{s\})\not=0\,.$}\par
\vskip 3mm
In other words, $T$ has a ``diagonal part''. Such operators necessarily 
have an atomic part, and so are isomorphisms on some subspace $L^1(A)$ 
([18], Theorem 5.5)$\,;$ this means that they are  Enflo operators ([8], 
Theorem 4.1). Translations are Enflo operators but not 
strong. The proof of Theorem 5.4 in [18] essentially shows$\,:$\par
\vskip 5mm
\noindent{\bf Proposition III.4.} {\it If $T$ is a strong Enflo operator, 
there exists a set $A$ of positive measure such that 
$$B_{L^1(A)}\ksubset \ind_A.T\big(B_{L^1(A)}\big)\ .$$}\par
\vskip 3mm
Indeed, if $T$ is a strong Enflo operator and 
$T^af(s)=\dis{\sum_{n=1}^{+\infty} a_n(s) f((\sigma_n(s))}$ 
represents the atomic part of $T,$ there exist an $n$ and a set $C$ of 
positive measure on which $a_n(s)\not=0$ and $\sigma_n(s)=s.$ It follows 
that, with the notation used in the proof of Theorem 5.4 of [18], 
$\sigma_n(B_{m,i}\cap B)\subseteq C_{m,i}\,,$ and this set 
$A=B_{m,i}\cap B$ gives the conclusion.\par
\vskip 5mm
\noindent{\bf Lemma III.5.} {\it Let $X$ be a subspace of $L^1(\Omega),$ 
and $T\colon L^1(A)\to L^1(\Omega)$ an operator, where 
$A\subseteq\Omega.$ Then $T\big(B_{L^1(A)}\big)\ksubset \ind_A.C_X$
if and only if there exists an operator 
$\widetilde T\colon L^1(A)\to X^{\perp\perp}$ such that 
$T=\ind_AP\widetilde T.$}\par
\vskip 3mm 
\noindent{\bf Proof.} The sufficiency comes directly from 
Bukhvalov-Lozanovski's theorem $P\big(B_{X^{\perp\perp}}\big)=C_X.$ To see 
the other direction, write $L^1(A)=\cup_{n\geq1} L_n$ with 
$L_n\subseteq L_{n+1}$ and $L_n$ isometric to $\ell_1^{\dim L_n}.$ Since 
$T\big(B_{L^1(A)}\big)\ksubset \ind_A.C_X$ there exist operators 
$T_n\colon L_n\to X$ such that $\Vert T_n\Vert\leq M$ and 
$d_m(\ind_AT_nf,Tf)\leq 2^{-n}$ for $f\in B_{L_n},$ where $d_m$ is a 
distance generating the topology $\tau_m.$ 
Taking a $w^\ast$-limit of the $T_n$'s along an ultrafilter gives 
$\widetilde T\colon L^1(A)\to X^{\perp\perp}$ such that 
$\ind_AP\widetilde T=T$. Indeed, if $F\in L^{1\ast\ast}$ is the 
{\it w}$^\ast$-limit of a bounded filter $(f_\alpha)_\alpha$ of functions 
in $L^1$, there exists a sequence $(c_k)_k$ of convex combinations of the 
$f_\alpha$'s such that $d_m(c_k,P(F))\dlowlim{\longrightarrow}
{\scriptscriptstyle{k\to+\infty}}0$ (see [16], Lemma IV.3.1).\hfill
$\qed$\par
\vskip 3mm
We can now state the following characterization of small subspaces.\par
\vskip 5mm
\noindent{\bf Proposition III.6.} {\it A subspace $X$ of $L^1$ is small if 
and only if no strong Enflo operator $T$ verifies 
$T\big( B_{L^1}\big)\ksubset C_X.$}
\par\vskip 3mm
\noindent{\bf Proof.} Suppose that $T$ is a strong Enflo operator such that 
$T\big( B_{L^1}\big)\ksubset C_X.$ Since by Proposition III.4, we have  
$B_{L^1(A)}\ksubset \ind_A.T\big(B_{L^1(A)}\big)$ for some set $A$ of 
positive measure, we obtain 
$B_{L^1(A)}\ksubset \ind_A.C_X\,,$ in contradiction with the smallness of 
$X.$\par
Suppose now that $X$ is not small, and let $A$ be of positive measure and 
such that $B_{L^1(A)}\ksubset \ind_A.C_X.$ By applying Lemma III.5 to the 
natural injection $j_A\colon L^1(A)\to L^1$ we obtain 
$\widetilde {j_A}\colon L^1(A)\to X^{\perp\perp}$ such that 
$j_A=\ind_AP\widetilde{j_A}.$ 
Defining $Tf=P\widetilde{j_A}(\ind_Af)$ for $f\in L^1,$ we get a strong 
Enflo operator $T\colon L^1\to L^1,$ since by uniqueness of the 
representation of $T$ we have $\mu_s^T=\delta_s$ for almost every 
$s\in A,$ and we have also
$T\big( B_{L^1(A)}\big)\ksubset C_X.$\hfill$\qed$\par
\vskip 5mm
Subspaces $X$ for which $C_X=B_X,$ i.e. nicely placed subspaces (see [10], 
D\'efinition 2), are of particular interest$\,:$ if $X$ is nicely placed and 
no operator from $L^1$ to $X$ is strongly Enflo, then $X$ is small. Since by 
Proposition III.4 the range of every strong Enflo operator contains an 
isomorphic copy of $L^1,$ we get$\,:$
\vskip 5mm
\noindent{\bf Proposition III.7.} {\it Every nicely placed subspace of $L^1$ 
which contains no subspace isomorphic to $L^1$ is small.}\par
\vskip 5mm
Note that the above condition is not necessary. For instance, the subspace 
of $L^1([0,1]^2)$ consisting of the functions which do not depend on the 
second variable is isometric to $L^1$ and nicely placed, but it is small. A 
special case of Proposition III.7 is$\,:$
\vskip 5mm
\noindent{\bf Proposition III.8.} {\it For every nicely placed Riesz subset 
$\Lambda$ of ${\ba Z}$, $L^1_\Lambda({\ba T})$ is small.}\par
\vskip 3mm
Recall that $\Lambda$ is a Riesz set if every measure whose Fourier 
transform is carried by $\Lambda$ is absolutely continuous, and it is nicely 
placed if $L^1_\Lambda$ is nicely placed in $L^1$ (see [11], Definition 1.4 
or [16], Definition IV.4.2). Proposition III.8 follows from the fact that 
$\Lambda$ is a Riesz set if and only if $L^1_\Lambda({\ba T})$ has 
Radon-Nikodym Property ([27]).\par
A stronger result will be shown below (Proposition III.10).\par
\vskip 5mm
Proposition III.6 gives also a stability result.\par
\vskip 5mm
\noindent{\bf Proposition III.9.} {\it Let $X$ and $Y$ be two small 
subspaces such that $X\cap Y=\{0\}$ and such that $X\oplus Y$ is closed. 
Then $X\oplus Y$ is small.}\par
\vskip 3mm
\noindent{\bf Proof.} Let $T\colon L^1\to L^1$ be such that 
$T\big(B_{L^1}\big)\ksubset C_{X\oplus Y}.$ Lemma III.5 with $A=\Omega$ 
gives 
$\widetilde T\colon L^1\to (X\oplus Y)^{\perp\perp}=X^{\perp\perp}\oplus 
Y^{\perp\perp}$ such that $T=P\widetilde T.$\par
Let $\pi\colon X^{\perp\perp}\oplus Y^{\perp\perp}\to X^{\perp\perp}$ be 
the projection. Writing $U=\pi\widetilde T$ and $V=(1-\pi)\widetilde T,$ 
we have $T=PU+PV,$ and $PU\big(B_{L^1}\big)\ksubset C_X$ and 
$PV\big(B_{L^1}\big)\ksubset C_Y.$ By Proposition III.6, $PU$ and $PV$ 
cannot be strong Enflo operators, so 
$\mu_s^{PU}(\{s\})=\mu_s^{PV}(\{s\})=0$ for 
almost all $s.$ Hence $\mu_s^T(\{s\})=0$ for almost all $s,$ and $T$ is 
not a strong Enflo operator. By Proposition III.6 again, we can conclude 
that $X\oplus Y$ is small.\hfill$\qed$\par 
\vskip 5mm
We now examine the translation invariant case. A subset $\Lambda$ of 
${\ba Z}$ is said to be {\sl semi-Riesz} ([37]) if every measure 
$\sigma$ whose Fourier transform is carried by $\Lambda$ has no atomic 
part. Every Riesz set is semi-Riesz, and Wiener's theorem ([20], p. 
42$\,;$ [13], p. 415)$\,:$
$$\lim_{N\to +\infty}{1\over 2N+1}\sum_{n=-N}^N
\vert \widehat\sigma(n)\vert^2=\sum_{t\in{\ba T}}\vert \sigma(\{t\}\vert^2$$
shows that every set $\Lambda$ with density zero, that is
$$d(\Lambda)=\dis{\lim_{N\to+\infty}{\sharp(\Lambda\cap \{-N,\cdots,N\})
\over 2N+1}}=0\,,$$ 
is semi-Riesz. The set 
$\Lambda=\big\{\dis{\,\sum_{k=1}^n}\varepsilon_k4^k\,;\ 
\varepsilon_k=-1,0,1, \ n\geq1\big\}$ is therefore semi-Riesz and not Riesz 
([37], p. 126). We have$\,:$\par
\vskip 5mm
\noindent{\bf Proposition III.10.} {\it Let $\Lambda$ be a nicely placed 
subset of ${\ba Z}\,.$ Then $L^1_\Lambda({\ba T})$ is small if and only if 
$\Lambda$ is a semi-Riesz set.}
\vskip 3mm
\noindent{\bf Proof.} Suppose that $L^1_\Lambda({\ba T})$ is not small, and 
let $T\colon L^1({\ba T})\to L^1({\ba T})$ be the operator constructed 
in the proof of Proposition III.6. We have readily 
$$\left\{\eqalign{&\mu_s^T(\{s\})=0\hskip 5mm{\rm if}\hskip 5mm 
s\not\in A\,;\cr
&\mu_s^T(\{s\})=1\hskip 5mm{\rm if}\hskip 5mm 
s\in A\,.}\right.$$\par
Let $\tau_t\colon L^1({\ba T})\to L^1({\ba T})$ be the translation given by
$$\tau_t f(u)=f(t+u)\ .$$\par
Since for every operator $U,$ we have 
$$\mu_s^{\tau_t^{-1}U\tau_t}(\{s\})=\mu_{s-t}^U(\{s-t\})\ ,$$
if we define $\widetilde T\colon L^1({\ba T})\to L^1({\ba T})$ by
$$\widetilde T=\int_{\ba T} \tau_t^{-1}T\tau_t\,dm(t)\ ,$$
we have, for all $s,$
$$\mu_s^{\widetilde T}(\{s\})=\int_{\ba T}\mu_{s-t}^T(\{s-t\})\,dm(t)\ ,$$
and so $\mu_s^{\widetilde T}(\{s\})=m(A)^{-1}.$\par
Since $\tau_t\widetilde T=\widetilde T\tau_t$ for every $t\in{\ba T},$ 
there exists a measure $\sigma$ such that 
$\widetilde T(f)=\sigma\ast f\,.$ We have then 
$\sigma(\{0\})=m(A)^{-1},$ so the atomic part of $\sigma$ is non zero, 
although $\sigma\in{\cal M}_\Lambda$.\par
Conversely, suppose that $L_\Lambda^1$ is small. Let $\sigma$ be a 
measure with spectrum in $\Lambda$ and with atomic part 
$$\sigma_a=\sum_{n=1}^{+\infty} a_n\delta_{x_n}\ .$$\par
The convolution operator defined by $Tf=\sigma\ast f\,$ maps $L^1({\ba T})$ 
into $L_\Lambda^1,$ so that $T\big(B_{L^1({\ba T})}\big)\ksubset 
B_{L_\Lambda^1}.$ On the other hand, its representing measure has the 
atomic part 
$$\mu_s^{T_a}=\sum_{n=1}^{+\infty} a_n\delta_{s-x_n}\ .$$
Since $L_\Lambda^1$ is small and nicely placed, $T$ cannot be a 
strong Enflo operator, so $\mu_s^T(\{s\})=0$ almost everywhere. This 
is possible only if $x_n\not=0$ or $a_n=0.$ Considering now, for each 
$p\geq1,$  the translated measure $\sigma\ast\delta_{-x_p}$ (which still 
have its spectrum in $\Lambda$) instead of $\sigma$ gives that $a_p=0,$ 
since its atomic part is  $\sum_{n=1}^{+\infty} a_n\delta_{x_n-x_p}\,.$ 
Hence $\sigma_a=0$ and we are done.\par\hfill$\qed$\par
\vskip 5mm
For every subset $\Lambda$ of ${\ba Z},$ denote by $[\Lambda]$ the smallest 
nicely placed subset of ${\ba Z}$ containing $\Lambda$ (see [11], p. 306). 
$[\Lambda]$ is contained in the Bohr-closure of $\Lambda$ ([11], Corollary 
2.6).\par
\vskip 3mm
\noindent{\bf Example.} Y. Meyer ([28], Th\'eor\`eme 6) showed that if 
$(n_k)_{k\geq1}$ grows fastly enough$\,:$\vskip -3mm
$$n_{k+1}>6(n_1+\cdots+n_k)\hskip 5mm {\rm and}\hskip 5mm 
\sum_{k=1}^{+\infty} {n_k\over n_{k+1}}<+\infty\,,$$ 
then the set 
$$\Lambda=\big\{\,\sum_{k=1}^K \varepsilon_k n_k\,;\ 
\varepsilon_k=-1,0,1,\ K\geq1\big\}$$ 
is closed for the Bohr topology$\,;$ hence it is nicely placed. Since it 
has density zero, $L^1_\Lambda({\ba T})$ is small, but $\Lambda$ is not a 
Riesz set, as the spectrum of a Riesz product. F. Parreau pointed out to us 
that the same conclusion holds for\vskip -2mm 
$$\Lambda=\big\{\,\sum_{k=1}^K \varepsilon_k 4^k\,;\ 
\varepsilon_k=-1,0,1,\ K\geq1\big\}\,.$$\par
We do not know whether $L^1_\Lambda$ contains a subspace isomorphic to 
$L^1$ whenever $\Lambda$ is the spectrum of a Riesz product. 
\vskip 5mm
\noindent{\bf Corollary III.11.} {\it If $L^1_\Lambda({\ba T})$ is not 
small, $[\Lambda]$ contains a translate of the spectrum of a Riesz product.}
\vskip 3mm
\noindent{\bf Proof.} If $L^1_\Lambda({\ba T})$ is not small, 
$L^1_{[\Lambda]}({\ba T})$ is {\it a fortiori} not small. Hence there is a 
measure $\sigma$ with spectrum in $[\Lambda]$ which has an atomic part, and 
Wiener's theorem shows that 
$\widehat\sigma(n)\mathop{{\hbox to 8mm{\rightarrowfill}}\llap
{\hbox{${\scriptstyle/}\hskip 3,5mm$}}}
\limits_{\scriptscriptstyle \vert n\vert\to+\infty} 0.$ 
Host-Parreau's theorem ([17]) ends the proof.\hfill
$\qed$\par
\vskip 1cm 




\def\N{{\rm I \kern-.25em N}}
\def\ind{{\rm 1\kern-.30em I}}
\def\dis{\displaystyle}
\def\wtd{\widetilde}
\def\dlowlim#1#2{\mathop{\rm #1}\limits_{#2}}

\magnification=1200

\noindent{\bf IV. LIFTING OF OPERATORS BETWEEN QUOTIENTS BY SMALL
SUBSPACES OF ${\bfmath L^1}$.}
\vskip 3mm 
In this section, we show that operators between quotients of $L^1$ by small 
subspaces lift to operators  on $L^1.$ We begin with a special case, which 
is much easier than the general one.\par
The following notion will be used$\,:$\par
\vskip 5mm
\noindent{\bf Definition IV.1.} {\it The operator $T\colon L^1\to L^1$ is 
said to be a {\sl Daugavet operator} if the Daugavet equation 
$\Vert I+\lambda T\Vert =1+\Vert T\Vert$ is fulfilled for every 
$\lambda$ with $\vert\lambda\vert=1.$}
\vskip 3mm
A. Plichko and M. Popov showed ([30], \S\kern.15em 9, Theorem 8) that every 
non-Enflo operator is Daugavet. In fact, more is true$\,:$ $T$ is a Daugavet 
operator whenever it is not strongly Enflo. Indeed, denoting the 
representing measures of $T$ by $\mu_s,$ $s\in[0,1],$ we have 
$$\Vert T\Vert=\sup_E{1\over m(E)}\int_0^1\vert\mu_s\vert(E)\,dm(s)$$
and
$$\Vert I+\lambda T\Vert=\sup_E{1\over m(E)}
\int_0^1\vert\delta_s+\lambda\mu_s\vert(E)\,dm(s)\,.$$ 
Since $T$ is not strongly Enflo, we have $\mu_s(\{s\})=0$ for almost all 
$s$, and the measures $\delta_s$ and $\mu_s$ are disjoint for these 
values. Hence, if for every $\varepsilon>0$, we choose a measurable set 
$E$ such that 
$$\Vert T\Vert\leq {1\over m(E)}\int_0^1\vert \mu_s\vert(E)\,dm(s)+
\varepsilon\,,$$
we get, since we may assume that $s\in E\,:$
$$1+\Vert T\Vert\leq {1\over m(E)}\int_0^1\vert \delta_s+\lambda
\mu_s\vert(E)\,dm(s)+\varepsilon\leq \Vert I+\lambda T\Vert+
\varepsilon\,,$$ 
which gives the result. D. Werner used a similar argument in ([36]).
\vskip 5mm
\noindent{\bf Proposition IV.2.} {\it Let $X$ and $Y$ be nicely placed 
subspaces of $L^1.$ Assume that $X$ and $Y$ are small, and that
$$\mathop{\rm dist}\big( L^1/X,L^1/Y\big)<2\ .$$\par
Then there exists an invertible operator $V\colon L^1\to L^1$ such that 
$V(X)=Y\,.$ 
Moreover, for $0\leq\alpha<1,$ if 
$$\mathop{\rm dist}\big( L^1/X,L^1/Y\big)=1+\alpha\,,$$
there is such a $V$ satifying $\Vert V\Vert.\Vert V^{-1}\Vert\leq 
\dis{1+\alpha\over 1-\alpha}\,\cdot$}\par
\vskip 3mm
\noindent{\bf Proof.} Let $S$ be an isomorphism from $L^1/X$ onto $L^1/Y$
such that $\Vert S\Vert.\Vert S^{-1}\Vert<2$. By the lifting property 
of nicely placed subspaces (see [23], Theorem 1, [16], Proposition 
IV.2.12, or [19], Proposition 2.1), there exist $U$ and $V$ such that
$$\left\{\eqalign{
Q_XU&=S^{-1}Q_Y\cr
Q_YV&=SQ_X}\right.$$
and such that $\Vert U\Vert=\Vert S^{-1}\Vert$ and $\Vert V\Vert=
\Vert S\Vert.$ Hence the following diagram is commutative$\,:$
$$\matrix{L^1&  
{\mathop{V}\limits_{\hbox to 1cm{\rightarrowfill}} 
\atop\raise 2mm\hbox{$
\dis{\mathop{\hbox to 1cm{\leftarrowfill}}\limits_{U}} $} } & L^1\cr
{\scriptstyle Q_X}\left\downarrow\vbox to 5mm{}\right. & & 
\left\downarrow\vbox to 5mm{}\right.  {\scriptstyle Q_Y}\cr
L^1/X & { \mathop{S}\limits_{\hbox to 1cm{\rightarrowfill}} \atop\raise 2mm
\hbox{$\dis\mathop{\hbox to 1cm{\leftarrowfill}}\limits_{S^{-1}}$} } & 
L^1/Y\cr}$$          
\par
For every $f\in L^1,$ we have
$$Q_YVUf=SQ_XUf=SS^{-1}Q_Yf=Q_Yf\ ,$$
hence
$$(I-VU)(L^1)\subseteq Y\ ,$$
and similarly $(I-UV)(L^1)\subseteq X.$\par
Since $X$ and $Y$ are nicely placed and small, the operators $I-UV$ and 
$I-VU$ are not strongly Enflo, and so are Daugavet operators. Therefore, 
we have 
$$\Vert I-UV\Vert=\Vert UV\Vert-1\hskip 0,5cm{\rm and}\hskip 0,5cm 
\Vert I-VU\Vert=\Vert VU\Vert-1\,.$$\par
Since $\Vert UV\Vert\leq\Vert U\Vert.\Vert V\Vert=\Vert S\Vert.
\Vert S^{-1}\Vert<2,$ it follows that $U$ and $V$ are invertible. 
Consequently $V(X)=Y$ since, for $f,g\in L^1$ such that $V(f)=g,$ we 
have
$$\eqalign{g\in Y\iff Q_Y(g)=0&\iff Q_YV(f)=0\cr
&\iff SQ_X(f)=0\iff Q_X(f)=0\iff f\in X\,.}$$\par
To conclude the proof, we observe that $\Vert V\Vert=\Vert S\Vert$ and 
$$V^{-1}=\Big[\sum_{k=0}^{+\infty}(I-UV)^k\Big]U\ ;$$
hence
$$\Vert V^{-1}\Vert \leq {\Vert S^{-1}\Vert\over 2-\Vert S\Vert.
\Vert S^{-1}\Vert}$$
since $\Vert U\Vert=\Vert S^{-1}\Vert$ and $\Vert I-UV\Vert\leq 
\Vert S\Vert.\Vert S^{-1}\Vert\,,$
and so $V$ works.\hfill$\qed$
\par\vskip 0,5cm
\noindent{\bf Remark.} Note that Proposition IV.2 applies in particular to 
nicely placed subspaces which do not contain $L^1$ (by Proposition III.7). 
If $L^1/X$ and $L^1/Y$ are isometric, we have $\alpha=0$ and, by the above, 
$V$ is actually an invertible isometry of $L^1.$\par 
\vskip 5mm
We now give the general statement.\par
\vskip 3mm
\noindent{\bf Theorem IV.3.} {\it Let $X$ and $Y$ be small subspaces 
of $L^1$ and suppose that there is an isomorphism $S\colon L^1/X\to L^1/Y$ 
such that $\Vert S\Vert$ and $\Vert S^{-1}\Vert$ are less than $1+\delta$ 
with $\delta<1/25.$\par
\ Then, there exists an invertible operator $U\colon L^1\to L^1$ such 
that\ \  
$\Vert U\Vert.\Vert U^{-1}\Vert\leq{(1+\delta)/(1-25\delta)}$ and 
$d_{\cal H}\big(U( B_X),B_Y\big)\leq 71\delta/(1-25\delta),$ 
where $d_{\cal H}$ denotes the Hausdorff distance.}
\vskip 3mm
\noindent{\bf Proof.} Denote by $Q_X$ and $Q_Y$ the projections onto the 
respective quotient spaces. We first construct a bidual diagram
$$\matrix{{L^1}^{\ast\ast}& \raise 2mm\hbox{$
{\mathop{W_1}\limits_{\hbox to 1cm{\rightarrowfill}} 
\atop {\dis{\mathop{\hbox to 1cm{\leftarrowfill}}\limits_{W_2} }} } $} 
& {L^1}^{\ast\ast}
\cr{\scriptstyle Q_X^{\ast\ast}}\left\downarrow\vbox to 5mm{}\right. & & 
\left\downarrow\vbox to 5mm{}\right.  {\scriptstyle Q_Y^{\ast\ast}}\cr
{L^1}^{\ast\ast}/X^{\perp\perp} 
&\raise 2mm\hbox{$ 
{\mathop{S^{\ast\ast}}\limits_{\hbox to 1cm{\rightarrowfill}} \atop
{\dis\mathop{\hbox to 1cm{\leftarrowfill}}\limits_{S^{-1\ast\ast}}} }$} 
& {L^1}^{\ast\ast}/Y^{\perp\perp}\cr}$$          
Write $L^1=\overline{\mathop{\bigcup}\limits_{n\geq1}L_n}$ where 
$\big(L_n\big)$ is an increasing sequence of finite-dimensional spaces 
isometric to $\ell_1^{\dim(L_n)}.$ For all $n\geq1,$ there exists 
$T_n\colon L_n\to L^1$ such that 
$$Q_Y T_n={SQ_X}_{\vert L_n}$$
and $\Vert T_n\Vert\leq (1+1/n)\,\Vert S\Vert\,.$
\vskip 5mm
\noindent{\bf Lemma IV.4.} {\it Let ${\cal U}$ be free ultrafilter on the 
integers. There exists a (w$^\ast$-w$^\ast$)-continuous operator 
$W_1\colon {L^1}^{\ast\ast}\to {L^1}^{\ast\ast}$ such that
$$\left\{\eqalign{&\Vert W_1\Vert\leq \Vert S\Vert\cr
&Q_Y^{\ast\ast}W_1=S^{\ast\ast}Q_X^{\ast\ast}\cr
&W_1(f)=F\ ,}\right.$$
for any 
$f\in\mathop{\bigcup}\limits_{n\geq1}L_n$ and  
$F={\it w}^\ast\hbox{\rm-}\dlowlim{\lim}{\cal U} T_n(f).$}
\vskip 3mm
\noindent{\bf Proof.} Define
${{\rm w}_1\colon \mathop{\bigcup}\limits_n L_n\to 
{L^1}^{\ast\ast}}$ 
by
$${\rm w}_1(f)={\it w}^\ast\hbox{\rm-}\dlowlim{\lim}{\cal U} T_n(f)\ .$$ 
\par
Clearly $\Vert{\rm w}_1\Vert\leq\Vert S\Vert$ and ${\rm w}_1$ can be 
extended to an operator $\overline{\rm w}_1\colon L^1\to{L^1}^{\ast\ast}.$ 
If now ${\Pi\colon (L^1)^{(4)}\to{L^1}^{\ast\ast}}$ denotes the 
canonical ({\it w}$^\ast$-{\it w}$^\ast$)-continuous projection
(``restriction to ${L^1}^\ast$''), then $W_1=\Pi\,
\overline{\rm w}_1^{\ast\ast}$ works. 
Indeed, for any $g\in L^1$ one has
$$Q_Y^{\ast\ast}\overline{\rm w}_1(g)=SQ_X(g)$$
and the equation $Q_Y^{\ast\ast}W_1=S^{\ast\ast}Q_X^{\ast\ast}$ follows by 
{\it w}$^\ast$-continuity of all relevant operators.\hfill$\qed$\par
\vskip 3mm
We get of course by the same token an operator $W_2$ with 
$Q_X^{\ast\ast}W_2={S^{-1}}^{\ast\ast}Q_Y^{\ast\ast}$ and similar 
properties.\par
\vskip 5mm
\noindent{\bf Lemma IV.5.} {\it Let $P\colon {L^1}^{\ast\ast}\to L^1$ be 
the natural projection. Then
$$\big[P(I-W_1W_2)\big]\big(B_{{L^1}^{\ast\ast}}\big)\subseteq 
\Vert I-W_1W_2\Vert\,\overline{B_Y}^{\tau_m}.$$}\par
\vskip 3mm
\noindent{\bf Proof.} For any $h\in {L^1}^{\ast\ast},$ we have
$$Q_Y^{\ast\ast}W_1W_2h=S^{\ast\ast}Q_X^{\ast\ast}W_2h=
S^{\ast\ast}{S^{-1}}^{\ast\ast}Q_Y^{\ast\ast}h=Q_Y^{\ast\ast}h\ ;$$
hence $Q_Y^{\ast\ast}(I-W_1W_2)=0$ and 
$$(I-W_1W_2)\big(B_{{L^1}^{\ast\ast}}\big)\subseteq \Vert I-W_1W_2\Vert
\,B_{Y^{\perp\perp}}.$$
The result follows by Bukhvalov-Lozanovski's theorem 
$P\big(B_{Y^{\perp\perp}}\big)=\overline{B_Y}^{\tau_m}.$\hfill$\qed$
\vskip 5mm
Then, writing $T=P(I-W_1W_2)_{\vert L^1},$ it follows from Lemma IV.5, the 
smallness of $Y$ and Proposition III.6 that $T$ is not a strong Enflo 
operator$\,;$ hence $-T$ satisfies the Daugavet equation, and we have
$$\eqalign{1+\Vert P(I-W_1W_2)_{\vert L^1}\Vert &=1+\Vert T\Vert 
=\Vert I-T\Vert\cr
&=\Vert (PW_1W_2)_{\vert L^1}\Vert\leq\Vert S\Vert.\Vert S^{-1}\Vert\leq 
1+\delta\ .}$$
Therefore
$$\Vert P(I-W_1W_2)_{\vert L^1}\Vert\leq \delta\ ,\leqno(1)$$
and hence, for any $f\in L^1,$ one has
$$\Vert P(I-W_1W_2)(f)\Vert_1\geq (1-\delta)\Vert f\Vert_1\ .$$
\par
On the other hand, one also has 
$$\Vert W_1W_2f\Vert_1\leq\Vert S\Vert.\Vert S^{-1}\Vert\,\Vert f\Vert_1
\leq (1+\delta)\Vert f\Vert_1\ .$$\par
Since $P$ is an L-projection, we have then
$$\Vert W_1W_2f- PW_1W_2f\Vert_1=
\Vert W_1W_2f\Vert_1-\Vert PW_1W_2f\Vert_1\leq 2\delta\Vert f\Vert_1\,.$$
\par
It follows from this and $(1)$ that 
$$\Vert(I-W_1W_2)_{\vert L^1}\Vert\leq 3\delta\,.$$\par
Since $W_1$ and $W_2$ are {\it w}$^\ast$-continuous, this gives
$$\Vert I-W_1W_2\Vert\leq 3\delta\,.\leqno(2)$$\par
An identical argument shows that
$$\Vert I-W_2W_1\Vert\leq 3\delta\,.\leqno(3)$$\par
\vskip 3mm
Since $\delta<1/3,$ $(2)$ and $(3)$ show that $W_1W_2$ and $W_2W_1$ are 
invertible. Hence $W_1$ and $W_2$ are themselves invertible.\par
\vskip 5mm
\noindent{\bf Lemma IV.6.} {\it In the notation of Lemma IV.4, one has
$$\mathop{\Vert\;\Vert\mathop{\hbox{\rm -}}{\rm diam}}{(L_f)}\leq 6\delta\,
\Vert W_2^{-1}\Vert\,\Vert f\Vert_1\,,$$
where $L_f$ is the set of all w$^\ast$-cluster points of the $T_nf$'s.}\par
\vskip 3mm
\noindent{\bf Proof.} Since $\Vert f-W_2W_1f\Vert_1\leq 
3\delta\,\Vert f\Vert_1,$ by $(2),$  one has
$$\Vert W_2^{-1}f-W_1f\Vert_1\leq 3\delta\,\Vert W_2^{-1}\Vert\,
\Vert f\Vert_1\,.$$
By Lemma IV.4, any $F\in L_f$ can be written $F=W_1(f)$ for an appropriate 
$W_1,$ and the conclusion follows.\hfill$\qed$
\vskip 5mm
Note that it follows from $(2)$ and $(3)$ that $\Vert W_2^{-1}\Vert
\leq \psi(\delta)=(1+\delta)/(1-3\delta).$\par
\vskip 5mm
\noindent{\bf Lemma IV.7.} {\it Let $Z$ be a Banach space such that 
$Z^{\ast\ast}=Z\oplus_1 Z_s$. Let $(z_n)_n$ be a sequence in 
$Z$ and $L$ be the set of all its w$^\ast$-cluster points. If
$$\mathop{\Vert\;\Vert\mathop{\hbox{\rm -}}{\rm diam}}{(L)}=\varepsilon\ ,$$
we have $\mathop{\rm dist}{(l,Z)}\leq\varepsilon$ for all $l\in L.$}
\vskip 3mm
\noindent{\bf Proof.} Let $l\in L.$ By 
translating if necessary, we may and do assume that $l\in Z_s.$\par
Let $K=\big(B_{Z^\ast},{\it w}^\ast\big).$\par
By ([6], Lemma III.2.4), we have 
$$\Vert l\Vert=\mathop{\rm dist}(l,Z)=\hat l=-\check l\ ,$$
where $\hat l$ and $\check l$ are respectively the {\it u.s.c.} and the 
{\it l.s.c.} hulls of $l$. Thus
$${\it w}^\ast\mathop{\hbox{\rm-}}\mathop{\rm Osc}(l_{\vert K})(y)
=2\Vert l\Vert\ .$$  
Pick $\alpha>0.$ For any $n\geq1,$ let
$$F_n=\{y\in K\,;\ \vert z_k(y)-z_j(y)\vert\leq \varepsilon+\alpha
\hskip 3mm \forall k,j\geq n\}\ .$$
It follows from the assumption that 
$$K=\bigcup_{n\geq1} F_n\ .$$
Hence, by Baire's lemma, there is a non-void {\it w}$^\ast$-open 
set $U\subseteq K$ and an $N\geq1$ such that $U\subseteq F_N.$ It 
follows that for all $y\in U,$ we have
$$\vert l(y)-z_N(y)\vert\leq \varepsilon+\alpha\ ,$$
and since $z_N$ is {\it w}$^\ast$-continuous, this implies that 
for all $y\in U$ we have 
$$\mathop{\rm Osc}(l_{\vert K})(y)\leq 2(\varepsilon+\alpha)\ ;$$
hence $\mathop{\rm dist} {(l,Z)}\leq \varepsilon+\alpha\,,$ which 
concludes the proof since $\alpha>0$ was arbitrary.\hfill$\qed$ 
\vskip 5mm
By Lemmas IV.6 and IV.7, we have
$$\Vert PF-F\Vert=\mathop{\rm dist}(F,L^1)
\leq 6\delta\psi(\delta)\,\Vert f\Vert_1$$
for every $F\in L_f\,.$\par
Once we know this, it follows that one has
$$\Vert (PW_i-W_i)_{\vert L^1}\Vert\leq 6\delta\psi(\delta)
\hskip 5mm(i=1,2)\ .\leqno(4)$$\par
We define now $U,V\colon L^1\to L^1$ by 
$$U=(PW_1)_{\vert L^1}\hskip 1cm V=(PW_2)_{\vert L^1}\ .$$\par
It follows from $(1),$ $(2),$ $(3)$ and $(4)$ that 
$$\Vert I-UV\Vert\leq {22\delta\over 1-3\delta}\hskip 1cm
\Vert I-VU\Vert\leq {22\delta\over 1-3\delta}\ {\raise 1mm\hbox{,}}$$
so $UV$ and $VU,$ and thus $U$ and $V$ as well, are invertible for 
$\delta<1/25.$ Moreover, $\Vert U^{-1}\Vert$ and $\Vert V^{-1}\Vert$ 
are less than $1/(1-25\delta).$\par
We now observe that since 
$$SQ_Xf=Q_Y^{\ast\ast}W_1f$$
for all $f\in L^1,$ one has $W_1(X)\subseteq Y^{\perp\perp}.$ It 
follows that $\mathop{\rm dist}{(Uf,Y^{\perp\perp})}\leq 
\Vert (PW_1-W_1)_{\vert L^1}\Vert,$ so by $(4)$
$$\mathop{\rm dist}{(Uf,Y)}\leq2\mathop{\rm dist}{(Uf,Y^{\perp\perp})}< 
{13\delta\over1-3\delta}$$
for any $f\in B_X.$ It follows that
$$U\big(B_X\big)\subseteq B_Y+{27\delta\over 1-3\delta}B_{L^1}.$$
\par{\openup 1mm
Since the same inclusion holds when we exchange $X$ and $Y$ and replace $U$ 
by $V,$ we can find, for every $g\in B_Y$, an $f\in X$ such that 
$\Vert Vg-f\Vert_1\leq {\raise -2pt\hbox{$\dis{13\delta\over 
1-3\delta}$}}\,$. 
Letting $h=f/\Vert f\Vert_1$, we obtain}
$$\Vert g-Uh\Vert_1\leq\Vert g-V^{-1}f\Vert_1+\Vert V^{-1}\Vert\,
\Vert I-UV\Vert\,\Vert f\Vert_1+\Vert U\Vert(\Vert f\Vert_1-1)$$
so
$$B_Y\subseteq U\big(B_X\big)+{71\delta\over 1-25\delta}B_{L^1}\,,$$ 
and this finishes the proof.\hfill$\qed$
\vskip 1cm
Let us say that $Y\subseteq L^1$ has the {\sl $\lambda$-lifting 
property} if for any $T\colon L^1\to L^1/Y,$ there exists 
$\widetilde T\colon L^1\to L^1$ such that $T=Q_Y\widetilde T$ and 
$\Vert \widetilde T\Vert\leq \lambda\Vert T\Vert.$ If $Y$ is 
complemented in $Y^{\ast\ast},$ it has the lifting property  (see 
[19], Lemma 2.1).\par
\vskip 5mm
\noindent{\bf Corollary IV.8.} {\it Let $X$ and $Y$ be small subspaces 
of $L^1$ such that there is an isomorphism ${S\colon L^1/X\to L^1/Y}$  
such that $\Vert S\Vert.\Vert S^{-1}\Vert\leq 1+\delta,$ and assume  
that $Y$ has the $\lambda$-lifting property.\par
Then for $\delta\leq \delta(\lambda),$ there is $T\colon L^1\to L^1$ 
invertible such that $T(X)=Y$ and 
$$\Vert T\Vert.\Vert T^{-1}\Vert\leq 1+K(\lambda)\delta$$
where $K(\lambda)$ only depends on $\lambda.$}\par
\vskip 3mm
\noindent{\bf Proof.} Since $SQ_X=Q_Y^{\ast\ast}{W_1}_{\vert L^1},$ by 
$(4)$ one has for $\delta\leq \delta_0\,:$
$$\Vert SQ_X-Q_YU\Vert_{L^1\to L^1/Y}\leq C\delta\ .$$\par
By assumption, there is $E\colon L^1\to L^1$ with 
$$\left\{\eqalign{&\Vert E\Vert\leq C\lambda\delta\cr
&Q_YE=SQ_X-Q_YU\ .}\right.$$\par
Now, if $T=U+E,$ one has $SQ_X=Q_YT,$ and the conclusion follows, since 
this implies, when $T$ is invertible, that $T(X)=Y.$\hfill$\qed$
\vskip 3mm
\noindent{\bf Remark IV.9.} There exist subspaces $X$ and $Y$ of infinite 
codimension in $c_0$ whose duals are isometric, but for which there is no 
invertible isometry of $\ell_1$ which maps $X^\perp$ onto $Y^\perp$. 
Indeed, every isometry from $\ell_1$ onto itself is 
{\it w}$^\ast$-continuous$\,;$ hence there should exist an isometry of 
$c_0$ mapping $X$ onto $Y$. Therefore, it suffices to construct infinite 
codimensional isometric subspaces $X$ and $Y$ of $c_0$ in such a way that $X$ 
contains an element supported by a single integer, and that the elements of 
$Y$ are all supported by an even number of integers. Since invertible 
isometries of $c_0$ leave the size of the supports invariant, such $X$ and 
$Y$ cannot be mapped into each other by an isometry of all $c_0$.\par  
So Proposition IV.2 with $\ell_1$ instead of $L^1$ is false in full 
generality. That shows that it is essential that the measure space be purely 
non-atomic, at least for defining a proper notion of ``small'' subspace.\par
\vskip 3mm
\noindent{\bf Remark IV.10.} It is asked in [19] whether, for two infinite 
Sidon sets $S_1$ and $S_2$, the spaces $L^1_{\wtd S_1}$ and 
$L^1_{\wtd S_2}$ are isomorphic. It would be true if one can lift 
isomorphisms from $L^1/L^1_{\wtd S_1}$ onto $L^1/L^1_{\wtd S_2}$. But 
$L^1_{\wtd S_1}$ and $L^1_{\wtd S_2}$ are ``large'' subspaces of $L^1$ and 
the present techniques are apparently not applicable to this question.\par 
\vskip 1cm





\def\dist{\mathop{{\rm dist}}\limits}
\def\dlowlim#1#2{\mathop{\rm #1}\limits_{#2}}
\def\esp{\mathop{\ba E}\nolimits}
\def\prob{\mathop{{\ba P}\kern 0pt}}
\def\dprob{\mathop{d_{\ba P}}}
\def\ind{{\rm 1\kern-.30em I}}
\def\dis{\displaystyle}
\def\N{{\rm I \kern-.20em N}}
\def\wtd{\widetilde}

\magnification=1200

\noindent{\bf V. CONSTRUCTION OF A PECULIAR SUBSPACE OF ${\bfmath L^1}$.}\par
\vskip 3mm
In Part II, we have seen a unified approach for Alspach's and Dor's 
theorems. In this section, we will see that the results of Part II have 
limitations. More precisely, we will construct a subspace $X$ of 
$L^1$ which will show that the nearness to isometries of the small 
into-isomorphisms from $L^1$ into $L^1$ (and of the small 
into-isomorphisms from finite dimensional subspaces of $L^1$ to $L^1$ 
([21])) does not keep on to take place for general subspaces of $L^1$. 
Moreover, the same space will show that extension of isometries from a 
subspace of $L^1$ to $L^1$ (see [15], [26], [31]) does not keep on 
either to take place with almost isometries. It shows also that Dor's 
Theorem B ([7]), which says that sequences of functions $f_n$ in $L^1$ 
which are equivalent to the canonical basis of $\ell_1$ are essentially 
supported by disjoint measurable sets cannot be improved by replacing the 
$f_n$'s by finite dimensional subspaces$\,:$ our space $X$ is, for every 
$\varepsilon>0$, $(1+\varepsilon)$-isomorphic to an $\ell_1$-sum of finite 
dimensional spaces, but $\Vert\esp^{\cal S}-Id\Vert_{{\cal L}(X,L^1)}\geq 
1$ for every $\sigma$-algebra ${\cal S}$ generated by a measurable 
partition. Note also that in our terminology, the space $X$ we construct is 
small and nicely placed.\par\vskip 1mm
We work in this section with real Banach spaces.
\vskip 3mm
\noindent{\bf Theorem V.1.} {\it There exists a subspace $X$ of $L^1=
L^1(\Omega,\Sigma,{\ba P})$ such that\par
\item{(a)} For every isometry $T\colon X\to L^1$ there is a unique isometry 
$\widetilde T\colon L^1\to L^1$ with $\wtd T_{\vert X}=T\,$;\par
\item{(b)} For every $\varepsilon>0$, there exists an isomorphism 
$T_\varepsilon\colon X\to L^1$ with $\Vert T_\varepsilon\Vert\,\Vert 
T_\varepsilon^{-1}\Vert\leq 1+\varepsilon$, but for which 
$\Vert T_\varepsilon-J\Vert \geq 1/2$ for every isometry $J\colon 
X\to L^1$.\par}
\vskip 3mm
Before producing the proof, we note that ({\it b}) and Alspach's theorem 
(Theorem II.7) show that the operators $T_\varepsilon$ do not extend to $L^1$ 
with a proper control on the norm of the extension.\par
\vskip 3mm
\noindent{\bf Proof.} The subspace $X$ is constructed in a very similar way 
than the one given in [12], Example 4.1.1, but we work with the stable 
variables themselves instead of working with their absolute value.\par
\vskip 3mm
The following lemma will be used\par
\vskip 5mm
\noindent{\bf Lemma V.2.} {\it Let $p\geq1$, $\varepsilon>0$ and 
$k,l\geq1$. There exists an $N\geq1$ such that for every subspace $F$ of 
$X=\ell_p^N$ with $\dim F\geq N-k$, there exist $l$ norm one vectors 
$y_1,\ldots,y_l$ of $\ell_p^N$ with disjoint support such that $\dist(y_j,F)
\leq\varepsilon$ $(1\leq j\leq l)$.}\par
\vskip 3mm
\noindent{\bf Proof.} Take $H$ bigger than $(1+2/\varepsilon)^k$ and $N=Hl$. 
Since $\dim {(X/Z)}=k$, indices $i_j$ ($Hh<i_{2h+1}<i_{2h+2}\leq H(h+1),\  
0\leq h\leq l$) can be found such that $\mathop{\rm dist}{(e_{2h+2}-
e_{2h+1},F)}\leq \varepsilon$, where $(e_n)_{1\leq n\leq N}$ is the 
canonical basis of $\ell_p^N$.\hfill$\qed$ 
\vskip 3mm
Recall first some probabilistic facts (see [3], pp. 60--61). A random variable 
$Z$ is called $p$-stable if its characteristic function is
$\widehat Z(t)\buildrel{\rm def}\over=\esp{(e^{itZ})}=
e^{-c_p\vert t\vert^p}$, 
where $c_p$ is the normalization constant ($\Vert Z\Vert_1=1$). Now, if 
$(Z_n)_{n\geq1}$ is a sequence of independent $p$-stable variables, its 
linear span $X_p$ is isometric to $\ell_p$, and $Z_1,\ldots,Z_n,\ldots$ 
correspond to the usual vector basis of $\ell_p$.\par
Since 
$c_p\dlowlim{\longrightarrow}{\scriptscriptstyle p\to1}0$,  
for every $\varepsilon>0$, there is an $\alpha>0$ such 
that if $1<p\leq \alpha\,:$
\vskip -2mm
$$\dprob{(f,0)}\leq \varepsilon\hskip 1cm\forall f\in B_{X_p}\ ,$$\par
\noindent 
where $\dprob$ defines the convergence in probability.\par
By the strong law of large numbers, we also have 
$\dis\dlowlim{\lim}{l\to+\infty}{1\over l}\sum_{j=1}^l\vert Z_j\vert=\ind$
almost surely.
\vskip 5mm
Now the space $X$ will be the closed linear span of the constant function 
$\ind$ and a sequence of independent $p_k$-stable variables with 
$p_k\dlowlim{\longrightarrow}{\scriptscriptstyle k\to+\infty}1$. More 
precisely, we choose inductively$\,:$\vskip 3mm 
\item {1)} a sequence $(p_k)_{k\geq1}$ converging to 1 and a sequence 
$(\varepsilon_k)_{k\geq1}$ converging to 0 such that
$$\dprob{(f,0)}\leq \varepsilon_k\hskip 1cm\forall f\in B_{X_{p_k}}\ ,$$ 
where $X_{p_k}$ is the linear span of independent $p_k$-stable variables$\,;$
these sequences are chosen so that the subspaces $X_1,X_2,\ldots$ 
constructed in the third step nearly create an $\ell_1$ sum. More 
explicitely, $\varepsilon_{k+1}$ is chosen so that the condition
$$\dprob{(f,0)}\leq \varepsilon_{k+1}\hskip 1cm
\forall f\in B_{X_{p_{k+1}}}$$
imply that for $g\in X_{k+1}$ and $f\in X_1+\cdots+ X_k$ we have\vskip -2mm 
$$\Vert f+g\Vert_1\geq (1-1/2^{k+2})\,(\Vert f\Vert_1+\Vert g\Vert_1)\,;$$
\par\vbox{}
\vskip-2mm
\item {2)} a sequence of integers $(l_k)_{k\geq1}$ such that
$\dis\Big\vert{1\over l_k}\sum_{j=1}^{l_k}\vert Y_j\vert-\ind\Big\vert
\leq 1/k$ almost surely for every independent $p_k$-stable variables 
$Y_1,\ldots,Y_{l_k}\,;$\par{\openup 1mm
\item {3)} a sequence of integers $(N_k)_{k\geq1}$ such that if $X_k$ is the 
linear span of the independent $p_k$-stable variables $Z_j, j\in I_k$ ($I_k
\subseteq\N, \vert I_k\vert=N_k$), we can find, by Lemma V.2, for every 
finite dimensional subspace $F$ of $X_k$, $l_k$ independent $p_k$-stable 
variables $Y_1,\ldots,Y_{l_k}$ in $X_k$ such that 
$\Vert\ \Vert_1{\rm -}\dist{(Y_j,F)}\leq 1/2^k$ for $1\leq j\leq l_k.$\par}
\vskip3mm{\openup 1mm
The sets $I_1,I_2,\ldots$ are chosen as successive intervals of $\N$ so we 
have $\N=\smash{\bigcup\limits_{\scriptscriptstyle k\geq1}}I_k.$ All 
the $p_k$-stable variables for the different values of $k$ are independently 
chosen.\par
Finally, we may, and do, suppose that the $\sigma$-algebra $\sigma(X)$ 
generated by $X$ is all $\Sigma$, since if not we replace 
$L^1(\Omega,\Sigma,{\ba P})$ by $L^1(\Omega,\sigma(X),{\ba P}_{\sigma(X)}).$ 
Since $\ind\in X$, Hardin's theorem ([15], Corollary 4.3) now gives part 
({\it a}) of the Theorem.\par}
\vskip 3mm
For the part ({\it b}), we state two facts. The first one follows directly 
from the construction and the almost isometric embedding of finite 
dimensional subspaces of $L^1$ into $\ell_1$.
\vskip 3mm
\noindent{\bf Fact 1.} {\it The Banach-Mazur distance from $X$ to the set of 
subspaces of $\ell_1$ is $1$}.\par
\vskip 3mm 
It follows that for every $\varepsilon>0$, there exist an operator 
$T_\varepsilon\colon X\to L^1$ such that $\Vert T_\varepsilon\Vert.
\Vert T_\varepsilon^{-1}\Vert\leq 1+\varepsilon$ and a $\sigma$-algebra 
${\cal S}_\varepsilon$ generated by disjoint measurable parts of $\Omega$ 
such that $\esp^{{\cal S}_\varepsilon}T_\varepsilon=T_\varepsilon$.\par
\vskip 5mm
\noindent{\bf Fact 2.} {\it For every sub-$\sigma$-algebra ${\cal S}$ 
of $\Sigma$ generated by disjoint measurable sets of $\Omega$, one has
$$\sup_{f\in B_X}\Vert \esp^{\cal S}f-f\Vert_1\geq1\ .$$}\par
\vskip 3mm
\noindent{\bf Proof.} Denote by $S_1,S_2,\ldots$ a measurable partition of
$\Omega$ generating ${\cal S}$, with $\prob{(S_i)}>0$ for all $i$. For 
every $k\geq1$, let 
$$F_k=\{f\in X_k\,;\ \int_{S_i}f\,d{\ba P}=0\hskip 5mm 
\forall i\leq k\}\ .$$\par
Let $\varepsilon>0.$ By the third step of the construction, since 
$\dim{(X_k/F_k)}\leq k$, and since 
$\inf\limits_{i\leq k} {\ba P}{(S_i)}>0$, there is an integer $k_0$ such 
that we can find, for every $k\geq k_0$, $l_k$ independent $p_k$-stable 
variables $Y_1,\ldots,Y_{l_k}$ in $X_k$ so that
$$\left\vert\,{1\over {\ba P}(S_i)}\int_{S_i}Y_j\,d{\ba P}\,\right\vert\leq 
\varepsilon\hskip 1cm\forall j\leq l_k,\ \forall i\leq k\,.\leqno(\ast)$$
\par
Let $U_k=S_1\cup\ldots\cup S_k$.\par
By $(\ast)$, one has
$$\Vert \esp^{\cal S}(Y_j)\,\ind_{U_k}\Vert_1\leq \varepsilon\hskip 1cm
\forall j\leq l_k\,.$$
But, by the step 2) of the construction, there is an integer $k_1$ such that 
$$\Big\Vert\, {1\over l_k}\sum_{j=1}^{l_k}\vert Y_j\vert\,\ind_{U_k}
\Big\Vert_1\geq 1-\varepsilon$$
for $k\geq k_1\,,$ and then there is a $j_k\leq l_k$ such that
$$\big\Vert\,\vert Y_{j_k}\vert\,\ind_{U_k}\big\Vert_1\geq 1-\varepsilon\,.$$
Hence
$$\Vert Y_{j_k}-\esp^{\cal S}Y_{j_k}\Vert_1\geq 
\Vert Y_{j_k}\,\ind_{U_k}-(\esp^{\cal S}Y_{j_k})\,\ind_{U_k}\Vert_1\geq 
1-2\varepsilon\,,$$
and this proves the Fact 2.\hfill$\qed$
\vskip 5mm
More generally, one has
\vskip 5mm
\noindent{\bf Lemma V.3.} {\it For every isometry $J\colon X\to L^1$ and 
every $\sigma$-algebra generated by a measurable partition, one has
$$\sup_{f\in B_X}\Vert \esp^{\cal S}Jf-Jf\Vert_1\geq1\ .$$}\par
\vskip 3mm
\noindent{\bf Proof.} From ({\it a}), one has an isometry $\widetilde J
\colon L^1\to L^1$ which extends $J$. There are measurable functions 
$a,\sigma$ on $\Omega$ (see Section II) such that
$$\widetilde Jf(s)=a(s)f(\sigma(s))\ .$$
Since $\Vert a\Vert_1=\Vert \widetilde J\ind\Vert_1=1$, the measure 
$$dQ=\vert a\vert\,d\/\prob{}$$
is a probability measure. Moreover,
$$\int_\Omega\vert 1+f(t)\vert\,d\/\prob(t)=
\int_\Omega\vert 1+f(\sigma(s))\vert\,dQ(s)\ ;$$
hence $f$ and $f(\sigma)$ have the same distribution ([15], Theorem 
1.1)$\,;$ so $f(\sigma)$ is $p$-stable in $L^1(Q)$ whenever $f$ is $p$-stable 
in $L^1(\prob)$.\par
Let ${\cal S}$ the $\sigma$-algebra generated by a measurable partition 
$(S_i)$ of $\Omega$ and let ${\cal S}'$ be a refinement of this partition, 
so that the sign of $a$ is constant on each $S_i'$ ($>0,$ $<0$, or $=0$). In 
particular $\mathop{{\rm supp}}{(a)}=\cup\{S_i'\,\ a\not\equiv0\  {\rm on}
\ S_i'\}$. 
\par
Let $\varepsilon'>0$. Reproducing the proof of the above Fact 2 in $L^1(Q)$, 
we obtain that if $S_i'$ ($i\leq k$) is disjoint from  
$\mathop{{\rm supp}}{(a)}$, 
there exist for every $k\geq k_0$, $Y_1,\ldots,Y_{l_k}\in X_k$ such that 
$$\left\vert\,{1\over Q(S_i')}\int_{S_i'}Y_j(\sigma)\,dQ\,\right\vert\leq 
\varepsilon'\hskip 1cm\forall i,\forall j\leq l_k\ ,$$
that is
$$\left\vert\,{1\over Q(S_i')}\int_{S_i'}Y_j(\sigma)\vert a\vert\,d\/\prob
\,\right\vert\leq 
\varepsilon'\hskip 1cm\forall i,\forall j\leq l_k\ .$$  
But $a$ has a constant sign on $S_i'$, so this writes
$$\left\vert\,{1\over Q(S_i')}\int_{S_i'}Y_j(\sigma)a\,d\/\prob\,\right\vert
\leq \varepsilon'\hskip 1cm\forall i,\forall j\leq l_k\ .$$
Therefore
$$\left\vert\,{1\over \prob{(S_i')}}\int_{S_i'}Y_j(\sigma)a\,d\/\prob\,\right
\vert\leq {Q(S_i')\over\prob{(S_i')}}\varepsilon'\hskip 1cm\forall i,
\forall j\leq l_k\ .$$
Consequently, for every $\varepsilon>0$, there is $k_0'$ such that for 
$k\geq k_0'$ there exist $Y_1,\ldots,Y_{l_k}\in X_k$ such that
$$\left\vert\,{1\over \prob{(S_i')}}\int_{S_i'}Y_j(\sigma)a\,d\/\prob\,\right
\vert\leq\varepsilon\hskip 1cm\forall i,\forall j\leq l_k\ .$$
Now, since indices $i\leq k$ for which $S_i'$ is disjoint from 
$\mathop{{\rm supp}}{(a)}$ do not matter, one has, as in the proof of Fact 
2, for $k$ with $U_k=S_1'\cup\ldots\cup S_k'\in{\cal S}\,:$
$$\Vert\esp^{\cal S}J(Y_j)\,\ind_{U_k}\Vert_1\leq
\Vert\esp^{\cal S'}J(Y_j)\,\ind_{U_k}\Vert_1\leq \varepsilon\hskip 1cm
\forall j\leq l_k\,.$$\par
On the other hand, 
$${1\over l_k}\sum_{j=1}^{l_k}\vert Y_j(\sigma)\vert\dlowlim{\longrightarrow}
{\scriptscriptstyle k\to+\infty}\ind\hskip 1cm\prob{\rm-a.s.}\,;$$
hence
$${1\over l_k}\sum_{j=1}^{l_k}\vert J(Y_j)\vert\dlowlim{\longrightarrow}
{\scriptscriptstyle k\to+\infty}\vert a\vert\hskip 1cm Q{\rm-a.s.}\ .$$\par
Since $\Vert a\Vert_{L^1(\prob)}=1$, it follows that
$$\mathop{\underline{\rm lim}}_{k\to+\infty}\Big\Vert\,{1\over l_k}
\sum_{j=1}^{l_k}
\vert J(Y_j)\vert\,\ind_{U_k}\Big\Vert_1\geq1\ .$$ 
There is therefore a $k_1$ such that, for every $k\geq k_1$, there is a 
$j_k\leq l_k$ such that
$$\Vert J(Y_{j_k})\,\ind_{U_k}\Vert_1\geq 1-\varepsilon\,,$$
and so
$$\Vert J(Y_{j_k})-\esp^{\cal S}J(Y_{j_k})\Vert_1\geq 
\Vert J(Y_{j_k})\,\ind_{U_k}-\esp^{\cal S}J(Y_{j_k})\,\ind_{U_k}\Vert_1
\geq 1-2\varepsilon\,,$$
and this ends the proof of the lemma.\hfill$\qed$
\vskip 5mm
We can now finish the proof of Theorem V.1.\par
\vskip 3mm
Let $\varepsilon>0$, $T_\varepsilon\colon X\to L^1$ with 
$\Vert T_\varepsilon\Vert.\Vert T_\varepsilon^{-1}\Vert\leq 1+\varepsilon$ 
and ${\cal S}_\varepsilon$ such that 
$\esp^{{\cal S}_\varepsilon}T_\varepsilon=T_\varepsilon$, as 
previously defined.\par
For every isometry $J\colon X\to L^1$, we have
$$\Vert\esp^{{\cal S}_\varepsilon}J-T_\varepsilon\Vert=\Vert 
\esp^{{\cal S}_\varepsilon}J-\esp^{{\cal S}_\varepsilon}T_\varepsilon\Vert
\leq \Vert J-T_\varepsilon\Vert\,,$$
so
$$\Vert \esp^{{\cal S}_\varepsilon}J-J\Vert\leq
\Vert \esp^{{\cal S}_\varepsilon}J-T_\varepsilon\Vert+
\Vert T_\varepsilon -J\Vert\leq 2\Vert T_\varepsilon-J\Vert\ .$$\par
But from Lemma V.3, $\Vert \esp^{{\cal S}_\varepsilon}J-J\Vert\geq 1$, 
and so $\Vert T_\varepsilon-J\Vert\geq 1/2.$\hfill$\qed$
\vskip 5mm
\noindent{\bf Remarks.} 1) There is a gap in the proof of the implication 
$(iii)\Rightarrow (iv)$ of Corollary 3.5 in [12]$\,;$ and we do not know 
whether one has $(i)\Rightarrow (iv)$ in that Corollary 3.5. Anyway 
Theorem V.1 shows that one cannot replace the distance $d_m$ in [12], 
Corollary 3.5 $(iv)$ by the $\Vert\ \Vert_1$ distance.\par
2) We mention here that it follows from L. Schwartz's thesis ([34]$\,;$ see 
[2], Theorem 4.2.5), that M\"untz spaces $M_1(\{t^{n_k}\})$ (with $\sum 
n_k^{-1}<+\infty$) are examples of subspaces of $L^1$ almost isometric to 
subspaces of $\ell_1$. The case $p>1$ has been noticed in [5] (see comments 
after Corollary 1.8).
\vskip 5mm
The derivation of Lemma V.3 from Lemma V.2 can actually be put in a more 
general frame. This is the content of the next proposition.
\vskip 5mm \goodbreak    
\noindent{\bf Proposition V.4.} {\it Let $X_1$ and $X_2$ be isometric 
subspaces of $L^1(\Omega_1,\Sigma_1,\mu_1)$ and \penalty -10000
$L^1(\Omega_2,\Sigma_2,\mu_2)$ respectively. Suppose that $X_1$ contains 
the constant functions, and that, for some $\varepsilon>0$, there is a 
$\sigma$-algebra ${\cal A}$ generated by a measurable partition of 
$\Omega_1$ for which
$$\Vert f-\esp^{\cal A}f\Vert_1\leq \varepsilon\Vert f\Vert_1
\hskip 1cm(\forall f\in X_1)\,.$$\par
Then for $\varepsilon'>\varepsilon$, there is a 
$(\Sigma_2\otimes {\cal B}{\rm or)-}$measurable partition 
of $\Omega_2\times[0,1]$ generating a $\sigma$-algebra ${\cal B}$ for 
which one has
$$\Vert g-\esp^{\cal B}g\Vert_1\leq \varepsilon'\Vert g\Vert_1
\hskip 1cm(\forall g\in X_2)\,.$$}\par
\vskip 3mm
In this statement, we identify $L^1(\Omega_2,\Sigma_2,\mu_2)$ to a 
subspace of 
$L^1(\Omega_2\times[0,1],\Sigma_2\otimes {\cal B}{\rm or},
\mu_2\otimes dt)=L^1(\wtd\Omega_2,\wtd\Sigma_2,\wtd\mu_2)$ 
by letting $\wtd g(\omega,t)=g(\omega)$ for $\omega\in\Omega_2$, 
$t\in[0,1]$ and $g\in L^1(\Omega_2,\Sigma_2,\mu_2)$.\par
\vskip 3mm
\noindent{\bf Proof.} Denote by $A_1,A_2,\ldots$ the partition generating 
${\cal A}$, and by $\sigma_1$ and $\sigma_2$ the sub-$\sigma$-algebras of 
$\Sigma_1$ and $\Sigma_2$ generated by $X_1$ and $X_2$ respectively. Let 
$U$ be an isometry from $X_1$ onto $X_2$. Since $\ind\in X_1$, Hardin's 
theorem ([15], Corollary 4.3) ensures the existence of an into isometry 
\vskip -1mm 
$$T=\wtd U\colon L^1(\sigma_1)\to L^1(\Sigma_2)$$
\vskip -1mm\noindent 
which extends $U$ and whose range is $L^1(\sigma_2)$.\par
$T$ can be written as $Tf=a(f\circ\tau)$ with $\tau$ an isomorphism between 
$(\Omega_2,\sigma_2)$ and $(\Omega_1,\sigma_1)$.\par
Write $Jf=f\circ\tau$. $J$ is positive, and so, since 
$\varphi_n=\esp^{\sigma_1}(\ind_{A_n})$ verifies\vskip -1mm 
$$\varphi_n\geq0\hskip 1cm {\rm and}\hskip 1cm 
\sum_{n\geq 1}\varphi_n=\ind\,,$$\vskip -1mm\noindent
we have 
$\psi_n=J\varphi_n\geq0$ and $\dis\sum_{n\geq1}\psi_n=\ind\,.$\par
Note that $\nu=\vert a\vert\mu_2$ is a probability measure and $J$ is an 
isometry from $L^1(\Omega_1,\sigma_1,\mu_1)$ to 
$L^1(\Omega_2,\sigma_2,\nu)$. We have$\,:$
\vskip 5mm
\noindent{\bf Lemma V.5.} {\openup 2mm{\it If $\psi_n$ are positive 
$\sigma_2$-measurable functions with $\sum_{n\geq1}\psi_n=\ind$, there 
is a measurable partition in sets $B_1,B_2,\ldots \in 
\sigma_2\otimes {\cal B}{\rm or}$ such that 
$\wtd\psi_n=\esp^{\sigma_2}(\ind_{B_n})$, for $n\geq1$.}\par}
\vskip 3mm
\noindent{\bf Proof.} Set
$$B_n=\big\{(\omega,t)\,;\ \sum_{k=1}^{n-1}\psi_k(\omega)\leq t<
\sum_{k=1}^n\psi_k(\omega)\,\big\}\ .$$ 
We have for every $B\in\sigma_2\,:$
$$\int_{\wtd B}\esp^{\sigma_2}\ind_{B_n}=\int_{\wtd B} \ind_{B_n}=
(\nu\otimes dt)(B_n\cap \wtd B)=\int_B m\big[B_n(\omega)\big]\,d\nu(\omega)
=\int_B\psi_n(\omega)\,d\nu(\omega)\,,$$
so $\wtd\psi_n=\esp^{\sigma_2}\ind_{B_n}.$\hfill$\qed$\par
\vskip 3mm
We may suppose now that $\wtd a$ has a constant sign on each $B_n$. Indeed, 
since $\wtd a$ is $\sigma_2$-measurable and $\tau$ is an isomorphism from 
$(\Omega_2,\sigma_2)$ onto $(\Omega_1,\sigma_1)$, we have, if 
$B_n^+=B_n\cap\{\wtd a\geq 0\}$ and $B_n^-=B_n\cap\{\wtd a<0\}$
$$\esp^{\sigma_2}(\ind_{B_n^+})=\ind_{\{\wtd a\geq 0\}}
\esp^{\sigma_2}(\ind_{B_n})=\ind_{\{\wtd a\geq0\}}\wtd\psi_n=
\ind_{\{\wtd a\geq0\}}(\wtd\varphi_n\circ\tau)={\wtd{\esp}}^{\sigma_1}
(\ind_{\{a\circ\tau^{-1}\geq0\}}\ind_{A_n})\circ\tau$$
so, we may cut each $A_n$ with the sets $\{a\circ\tau^{-1}\geq 0\}$ 
and $\{a\circ\tau^{-1}<0\}$ and the $B_n^+$'s and the $B_n^-$'s will 
correspond to the $A_n^+$'s and the $A_n^-$'s.\par 
Now, letting 
$$a_n(f)={1\over \mu_1(A_n)}\int_{A_n}f\,d\mu_1\ ,$$
one has
$$\leqalignno{\mu_1(A_n)&=\int_{\Omega_1}\ind_{A_n}\,d\mu_1
=\int_{\Omega_1}\varphi_n\,d\mu_1=\Vert \varphi_n\Vert_1
=\Vert \psi_n\Vert_{L^1(\nu)} &1)\cr
&=\int_{\Omega_2}\psi_n\,d\nu=
\int_{\wtd\Omega_2}\ind_{B_n}\,d\wtd\nu=\wtd\nu(B_n)\ ,}$$
and 
$$\leqalignno{\int_{A_n}f\,d\mu_1
&=\int_{\Omega_1}f\ind_{A_n}\,d\mu_1=
\int_{\Omega_1}f\varphi_n\,d\mu_1=
\int_{\Omega_2}J(f\varphi_n)\,d\nu& 2)\cr
&\hskip 4cm\hbox{since $Jf=f\circ\tau$ and $\mu_1=\tau^\ast(\nu)$}\cr
&=\int_{\Omega_2}(Jf)\psi_n\,d\nu=\int_{B_n}\wtd Jf\,d\wtd\nu\ .\cr}$$ 
Hence
$$a_n(f)={1\over \wtd\nu(B_n)}\int_{B_n}\wtd Jf\,d\wtd\nu\ ,$$
and then, denoting by ${\cal B}$ the $\sigma$-algebra generated by 
$B_1,B_2,\ldots\,:$
$$\eqalign{\Vert \wtd Jf-\esp^{\cal B}(\wtd Jf)\Vert_{L^1(\tilde\nu)}
&=\sum_{n\geq1}\int_{B_n}\vert \wtd Jf-a_n(f)\ind\vert\,d\wtd\nu=
\sum_{n\geq1}\int_{B_n}\vert \wtd f\circ\tau-a_n(f)\ind)\vert\,d\wtd\nu\cr
&=\sum_{n\geq1}\int_{\Omega_2} \vert f\circ\tau-a_n(f)\ind\vert\psi_n\,d\nu
\cr
&=\sum_{n\geq1}\int_{\Omega_2} \vert f\circ\tau-a_n(f)\ind\vert
(\varphi_n\circ\tau)\,d\nu\cr
&=\sum_{n\geq1}\int_{\Omega_1} \vert f-a_n(f)\ind\vert\varphi_n\,d\mu_1
=\sum_{n\geq1}\int_{A_n} \vert f-a_n(f)\ind\vert\,d\mu_1\cr
&=\Vert f-\esp^{\cal A}f\Vert_1\leq\varepsilon\Vert f\Vert_1\ .}$$\par
\vskip 3mm
Now, going back to $T$, we refine ${\cal B}$ in such a way that 
$\Vert \wtd a-\esp^{\cal B}\wtd a\Vert_{L^1(\tilde\mu_2)}
\leq \alpha=\varepsilon'-\varepsilon$. We keep the inequality 
$\Vert \wtd Jf-\esp^{\cal B}(\wtd Jf)\Vert_{L^1(\tilde\nu)}\leq\varepsilon
\Vert f\Vert_1$. Moreover, this refinement can be made by taking the 
intersection with $\sigma_2$-measurable sets, so, as above, the new $B_n$'s 
will correspond to an appropriate refinement of the $A_n$'s.\par
Set $\varepsilon_n=1$ if $a\geq0$ on $B_n$ and $\varepsilon_n=-1$ if 
$a<0$ on $B_n$. One has
$$\eqalign{\int_{B_n}\wtd a\,d\wtd\mu_2
&=\varepsilon_n\int_{\Omega_2}\vert a\vert\psi_n\,d\mu_2
=\varepsilon_n\int\vert a\vert(\varphi_n\circ\tau)\,d\mu_2
=\varepsilon_n\int_{\Omega_2}\vert T\varphi_n\vert\,d\mu_2\cr
&=\varepsilon_n\Vert T\varphi_n\Vert_{L^1(\mu_2)}=\varepsilon_n
\Vert \varphi_n\Vert_{L^1(\mu_1)}=\varepsilon_n\mu_1(A_n)\,,}$$
and so\ \ $\esp^{\cal B}\wtd a={\dis \sum_{n\geq1}\varepsilon_n
{\mu_1(A_n)\over \wtd\mu_2(B_n)}\ind_{B_n}}.$\par
Now
$$\eqalignno{
\Vert \wtd Tf-\esp^{\cal B}\wtd Tf\Vert_{L^1(\tilde\mu_2)}
&=\int_{\wtd\Omega_2}\Big\vert \wtd a(\wtd f\circ\tau)-
\sum_{n\geq1}\Big({1\over \wtd\mu_2(B_n)}\int_{B_n}\wtd a(\wtd f\circ\tau)
\Big)\ind_{B_n}\Big\vert\,d\wtd\mu_2\cr
&\leq \int_{\wtd\Omega_2}\Big\vert \wtd a(\wtd f\circ\tau)-
\sum_{n\geq1}\Big({\wtd a\over \mu_1(A_n)}\int_{B_n}\vert\wtd a\vert
(\wtd f\circ\tau)_,d\wtd\mu_2\Big)\ind_{B_n}\Big\vert\,d\wtd\mu_2\cr
&\hskip 1cm+\int_{\wtd\Omega_2}\Big\vert \sum_{n\geq1}
\Big({\wtd a\over \mu_1(A_n)}
\int_{B_n}\vert\wtd a\vert(\wtd f\circ\tau)\,d\wtd\mu_2\Big)\ind_{B_n}\cr
&\hskip 4cm-\sum_{n\geq1}\Big({1\over \wtd\mu_2(B_n)}
\int_{B_n}\wtd a(\wtd f\circ\tau)
\,d\wtd\mu_2\Big)\ind_{B_n}\Big\vert\,d\wtd\mu_2\cr
&\leq\int_{\wtd\Omega_2}\Big\vert \wtd f\circ\tau-
\sum_{n\geq1}\Big({1\over \wtd\nu(B_n)}\int_{B_n}(\wtd f\circ\tau)
\,d\wtd \nu\Big)\ind_{B_n}\Big\vert\,d\wtd\nu\cr
&\hskip 7mm+\int_{\wtd\Omega_2}\sum_{n\geq1}\Big({1\over\mu_1(A_n)}
\int_{B_n}\vert \wtd a\vert(\wtd f\circ\tau)\,d\wtd\mu_2\Big)
\Big(\wtd a-\varepsilon_n{\mu_1(A_n)\over \wtd\mu_2(B_n)}\Big)\ind_{B_n}
\,d\wtd\mu_2\cr
&\leq\Vert \wtd Jf-\esp^{\cal B}\wtd Jf\Vert_{L^1(\wtd\nu)}+
\int_{\wtd\Omega_2}\sum_{n\geq1}\Vert f\Vert_1
\Big\vert \wtd a-\varepsilon_n{\mu_1(A_n)\over\wtd\mu_2(B_n)}\Big\vert
\ind_{B_n}\,d\wtd\mu_2\cr
&\leq \varepsilon\Vert f\Vert_1+\Vert f\Vert_1\sum_{n\geq1}\int_{B_n}
\Big\vert\wtd a-\varepsilon_n{\mu_1(A_n)\over\wtd\mu_2(B_n)}\Big\vert
\,d\wtd\mu_2\cr
&\leq \varepsilon\Vert f\Vert_1+\Vert f\Vert_1\Vert\wtd a-\esp^{\cal B}
\wtd a\Vert_{L^1(\wtd\mu_2)}\cr
&\leq (\varepsilon+\alpha)\Vert f\Vert_1=\varepsilon'\Vert f\Vert_1\ .
& \qed}$$ 

\vskip 5mm
\noindent{\bf Remarks.} 1) A slight variation in Proposition V.4 is$\,:$ set 
$C_n=\tau^{-1}(A_n)$. The $C_n$'s are in general not in $\Sigma_2$, but if 
${\cal C}$ is the $\sigma$-algebra generated by them, the map 
$\tau\colon\Omega_2\to\Omega_1$ is 
($\sigma_2\vee{\cal C}$-$\sigma_1\vee{\cal A}$)-bi-measurable (see in 
[14] a description of these $\sigma$-algebras), and 
$\widetilde\mu_2(B)=\mu_1[\tau(B)]$ defines a measure on 
$(\Omega_2,\sigma_2\vee{\cal C})$, whose restriction to $\sigma_2$ is 
equal to the restriction of $\mu_2$ to $\sigma_2$. One has then an 
isometry $\widetilde J\colon L^1(\Omega_1,\sigma_1\vee{\cal A},\mu_1)
\to L^1(\Omega_2,\sigma_2\vee {\cal C},\widetilde\mu_2)$ defined by 
$\widetilde J(f)=f\circ \tau$ such that $\widetilde JX_1=JX_1$. It 
follows that 
$\Vert Jf-\esp^{\cal C}Jf\Vert_1\leq\varepsilon\Vert f\Vert_1\,.$\par
It would be interesting to find an intrinsic characterization of these 
spaces. This problem may be connected to the fact that though there 
are sufficient conditions to have uniform convergence of martingales (see 
[29]), it seems there is at the present time no available necessary 
condition.\par
2) Since the space $X$ and $T_\varepsilon(X)$ are small nicely placed 
subspaces (by Proposition III.7), it follows from Theorem V.1 and 
Proposition IV.2 that there is $\alpha>0$ such that the quotient spaces 
$L^1/X$ and $L^1/T_\varepsilon(X)$ have Banach-Mazur distance greater than 
$(1+\alpha)$ for all $\varepsilon>0\,.$ This implies of course that the 
operators $T_\varepsilon$ do not extend to isomorphisms $U$ of $L^1$ such 
that the norms of $U$ and its inverse are small.\par
\vfill\eject



\magnification=1200
\noindent{\bf REFERENCES.}\par
\vskip 5mm
\noindent \item {1.} D.E. Alspach, {\it Small into isomorphisms on $L_p$ 
spaces}, Illinois J. Math 27 (1983) 300--314.\par
\smallskip
\noindent \item {2.} P. Borwein and T. Erd\'elyi, {\it Polynomials and 
polynomial inequalities}, Graduate Texts in Math. 161 (1995).
\smallskip
\noindent \item {3.} J. Bourgain and H.P. Rosenthal, {\it Martingales
valued in certain subspaces of $L^1$}, Isra\"el J. Math. 37 (1980) 54--75
\smallskip
\noindent \item {4.} A.V. Bukhvalov and G. Lozanovskii, {\it On sets closed 
in measure}, Trans. Moscow Math. Soc. 2 (1978) 127--148.
\smallskip
\noindent \item {5.} F. Delbaen, H. Jarchow and A. 
Pel\llap{\raise 1mm\vbox{\hbox{\fiverm/}}}czy\'nki, {\it Subspaces of $L_p$ 
Isometric to Subspaces of $\ell_p$}, Positivity ({\it to appear}).
\smallskip
\noindent \item {6.} R. Deville, G. Godefroy and V. Zizler, {\it Smoothness 
and Renormings in Banach Spaces}, Pitman Monographs in Pure and Applied 
Mathematics 64 (1993).
\smallskip
\noindent \item {7.} L.E. Dor, {\it On projections in $L^1$}, Annals of 
Math. 102 (1975) 463--474.\par
\smallskip
\noindent \item {8.} P. Enflo and T.W. Starbird, {\it Subspaces of $L^1$ 
containing $L^1$}, Studia Math. 65 (1979) 203--225.\par
\smallskip
\noindent \item {9.} H. Fakhoury, {\it Repr\'esentations d'op\'erateurs 
\`a valeurs dans $L^1(X,\Sigma,\mu)$}, Math. Ann. 240 (1979) 203--212.\par 
\smallskip
\noindent \item {10.} G. Godefroy, {\it Sous-espaces bien dispos\'es de 
$L^1\ -$ Applications}, Trans. Amer. Math. Soc. 286 (1984) 227--249.
\smallskip
\noindent \item {11.} G. Godefroy, {\it On Riesz subsets of abelian discrete 
groups}, Isra\"el J. Math. 61 (1988) 301--331.
\smallskip
\noindent \item {12.} G. Godefroy, N.J. Kalton and D. Li, {\it On subspaces 
of $L^1$ which embed into $\ell_1$}, Journal f\"ur die reine und angew. 
Math. 471 (1996) 43--75.\par
\smallskip
\noindent \item {13.} C.C. Graham and O.C. McGehee, {\it Essays 
in Commutative Harmonic Analysis}, \penalty -10000Springer-Verlag (1979).\par
\smallskip
\noindent \item {14.} A. Hanen and J. Neveu, {\it Atomes conditionnels d'un 
espace de probabilit\'e}, Acta Math. Acad. Sci. Hung. 17 (1966) 443--449.
\smallskip
\noindent \item {15.} C.D. Hardin Jr., {\it Isometries on Subspaces of 
$L^p$}, Indiana Univ. Math. J. 30 (1981) 449--465.\par
\smallskip
\noindent \item {16.} P. Harmand, D. Werner and W. Werner, {\it M-ideals in 
Banach spaces and Banach algebras}, Lecture Notes in Math. 1547 (1993), 
Springer-Verlag.
\smallskip
\noindent \item {17.} B. Host and F. Parreau, {\it Sur les mesures dont la 
transform\'ee de Fourier ne tend pas vers z\'ero \`a l'infini}, Colloq. 
Math. 41 (1979) 285--289.
\smallskip
\noindent \item {18.} N.J. Kalton, {\it The Endomorphisms of $L_p$ 
$(0\leq p\leq 1)$}, Indiana Univ. Math. J. 27 (1978) 353--381.\par
\smallskip
\noindent \item {19.} N.J. Kalton and A. 
Pel\llap{\raise 1mm\vbox{\hbox{\fiverm/}}}czy\'nski, {\it Kernels 
of surjections from ${\cal L}_1$-spaces with an application to Sidon 
sets}, Math. Ann. 309 (1997) 135--158.
\smallskip
\noindent \item {20.} Y. Katznelson, {\it An introduction to Harmonic 
Analysis}, Wiley (1976).
\smallskip
\noindent \item {21.} A.L. Koldobskii, {\it Convolution metric in the space 
of measures and $\varepsilon$-isometries in $L_p$}, J. Soviet Math. 44 
(1989) 852--855.\par
\smallskip
\noindent \item {22.} E. Lacey, {\it The Isometric Theory of Classical 
Banach Spaces}, Springer-Verlag (1974). 
\smallskip
\noindent \item {23.} D. Li, {\it Lifting Properties for Some Quotients of 
$L^1$-Spaces and Other Spaces \penalty-10000 $L$-summand in Their Bidual}, 
Math. Z. 199 (1988) 321--329.\par
\smallskip
\noindent \item {24.} J. Lindenstrauss and H.P. Rosenthal, {\it 
Automorphisms in $c_0$, $\ell_1$ and $m$}, Isra\"el J. Math. 7 (1969) 
227--239. 
\smallskip
\noindent \item {25.} J. Lindenstrauss and L. Tzafriri, {\it Classical 
Banach Spaces I and II}, Classics in Math., Springer (1996).
\smallskip
\noindent \item {26.} W. Lusky, {\it Some Consequences of Rudin's Paper 
``$L_p$-Isometries and Equimeasurability''}, Indiana Univ. Math. J. 27 
(1978) 859--866.\par
\smallskip
\noindent \item {27.} F. Lust-Piquard, {\it Ensembles de Rosenthal et 
ensembles de Riesz}, Comptes Rendus Acad. Sciences (Paris) S\'erie A 282 
(1976) 833--835.
\smallskip
\noindent \item {28.} Y. Meyer, {\it Spectres des mesures et mesures 
absolument continues}, Studia Math. 30 (1968) 87--99.
\smallskip
\noindent \item {29.} J. Neveu, {\it Note on the tightness of the metric on 
the set of complete sub $\sigma$-algebras on a probability space}, Ann. 
Math. Statistics 43 (1972) 1369--1371.
\smallskip
\noindent \item {30.} A.M. Plichko and M.M. Popov, {\it Symmetric function 
spaces on atomless probability spaces}, Dissertationes Math. 306 (1990).
\smallskip
\noindent \item {31.} A.I. Plotkin, {\it Continuation of $L^p$ isometries}, 
J. Soviet Math. 2 (1974) 143--165.
\smallskip
\noindent \item {32.} H.P. Rosenthal, {\it Embeddings of $L^1$ in $L^1$}, 
Contemp. Math. 26 (1984) 335--349.
\smallskip
\noindent \item {33.} G. Schechtman, {\it Almost isometric $L_p$ subspaces 
of $L_p(0,1)$}, Journal London Math. Soc. (2) 20 (1979) 516--528.
\smallskip
\noindent \item {34.} L. Schwartz, {\it \'Etude des sommes 
d'exponentielles}, Hermann (1943).
\smallskip
\noindent \item {35.} M. Talagrand, {The three-space problem for $L^1$}, J. 
Amer. Math. Soc. 3 (1990) 9--29.
\smallskip
\noindent \item {36.} D. Werner, {\it An elementary apprroach to the 
Daugavet equation}, {\it in ``Interaction between Functional Analysis, 
Harmonic Analysis, and Probability''} (N.J. Kalton, E. Saab, and S. 
Montgomery-Smith, Eds), Proc. Conf. Columbia 1994, Lecture Notes in Pure 
and Applied Mathematics, Vol. 175, pp. 449--454, M. Dekker (1996).
\smallskip
\noindent \item {37.} D. Werner, {\it The Daugavet Equation for Operators on 
Function Spaces}, J. Funct. Anal. 143 (1997) 117--128.\par
\vskip 1cm
\noindent {\it G. Godefroy}, Universit\'e Pierre et Marie Curie, Equipe 
d'Analyse, Tour 46-0, Case 186, 4 place Jussieu, F-75252 Paris Cedex 05, 
e-mail gig@ccr.jussieu.fr
\medskip
\noindent {\it N.J. Kalton}, University of Missouri-Columbia, Department of 
Mathematics, Columbia, MO 65211 USA, e-mail nigel@math.missouri.edu
\medskip
\noindent {\it D. Li}, Universit\'e d'Artois, P\^ole de Lens, Facult\'e des 
Sciences Jean Perrin, Rue Jean Souvraz, SP 18, F-62307 Lens Cedex, e-mail 
daniel.li@euler.univ-artois.fr
\par\vfill\eject


\bye